\documentclass[2p]{article}
\usepackage{amssymb,amsmath,amsthm}
\usepackage{indentfirst}
\usepackage{exscale}
\usepackage{relsize}
\usepackage[numbers,sort&compress]{natbib}

\usepackage{geometry}
\usepackage{color}

\theoremstyle{definition}

\allowdisplaybreaks
\theoremstyle{remark}

\numberwithin{equation}{section}
 \UseRawInputEncoding

\begin{document}

\title{\Large\bf{ Multiple solutions for a class of nonhomogeneous elliptic systems with Dirichlet boundary or Neumann boundary }
 }
\date{}
\author {{Xiaoli Yu}$^{1}$,{\ Xingyong Zhang} $^{1,2}$\footnote{Corresponding author, E-mail address: zhangxingyong1@163.com}  \\
{\footnotesize $^{1}$Faculty of Science, Kunming University of Science and Technology, Kunming, Yunnan, 650500, P.R. China.}\\
       {\footnotesize Kunming, Yunnan, 650500, P.R. China.}\\
       {\footnotesize $^{2}$Research Center for Mathematics and Interdisciplinary Sciences, Kunming University of Science and Technology,}\\
       {\footnotesize Kunming, Yunnan, 650500, P.R. China.}\\
 }
 \date{}
 \maketitle

 \begin{center}
 \begin{minipage}{15cm}

\small  {\bf Abstract:}
In this paper, we mainly establish the existence of at least three non-trivial solutions for a class of nonhomogeneous quasilinear elliptic systems with Dirichlet boundary  value or Neumann boundary value  in a bounded domain $\Omega\subset\mathbb{R}^N $ and $N\geq 1$.  We exploit  the method which  is based on \cite{Bonanno Bisci 2011}. This method let us  obtain the concrete open interval about the parameter $\lambda$. Since  the quasilinear term depends on $u$ and $\nabla u$, it is necessary for our proofs to use the theory of  monotone operators and the skill of adding one dimension to space.

 \par
 {\bf Keywords:}
three non-trivial solutions; concrete open interval; Dirichlet boundary  value;  Neumann boundary value; quasilinear scalar field systems.
\par
 {\bf 2020 Mathematics Subject Classification.} 35J20; 35J50; 35J62.
\end{minipage}
 \end{center}
  \allowdisplaybreaks
\section{Introduction and main results }
\label{section 1}
 \vskip2mm
 \noindent
 Stuart \cite{C.A. Stuart}  established the result of  at least two non-negative weak solutions for the following  nonhomogeneous problem:
\begin{eqnarray}
\label{af2}
 \begin{cases}
 -\text{div}\left(\phi\left(\dfrac{u^2+|\nabla u|^2}{2}\right)\nabla u\right)+\phi\left(\dfrac{u^2+|\nabla u|^2}{2}\right)u
=\lambda u+h,\;\;\;\;\;\; x\in \;\Omega ,\\
  u(x)= 0,\;\;\;\;\;\; x\in \;\partial\Omega,
   \end{cases}
\end{eqnarray}
where  $h\in L^2(\Omega),$ $h\geq 0$, $\Omega$ is a bounded domain of $\mathbb{R}^N,$ $N\geq 1$,  $\lambda \in \mathbb{R}$, $\phi\in C([0,+\infty),\mathbb{R})$, $\Phi(s)=\int_0^s\phi(\tau)d\tau$ and satisfies\\
$(\Phi_1)$  $\Phi(0)=0$, $\phi=\Phi'$ is non-increasing on $[0,+\infty)$ and $\phi(\infty)=\lim_{s\rightarrow\infty}\phi(s)>0;$ \\
$(\Phi_2)$ putting $\Phi(s^2)=g(s),$ there exists $\rho>0$ such that\\
$$ g(t)\geq g(s)+g'(s)(t-s)+\rho(t-s)^2,\;\;\;\;\mbox{for all}\; s,t\geq 0;$$\\
$(\Phi_3)$ $\lim_{s\rightarrow\infty} (\Phi(s)-\Phi'(s)s)<+\infty$.\\
$(\Phi2)$ is used to ensure the ellipticity of $(\ref{af2}).$  More specifically, when $\phi(\infty)+\phi(\infty)\lambda_1<\lambda<\phi(0)+\phi(\infty)\lambda_1$, he obtained  results that  the equation $(\ref{af2})$ has a trivial solution and another non-negative non-trivial weak solution for $h\equiv0$ and if   there exists $H_\lambda$ and $0<\|h\|_2=[\int_\Omega h^2dx]^{1/2}<H_\lambda $,   the equation $(\ref{af2})$ has at least two distinct  non-negative weak solutions. He claimed that  in both cases one is a strict local minimum solution and the other is  a mountain pass type solution.
\par
In the equation (\ref{af2}), the quasilinear term $\phi$ depends on $\dfrac{u^2+|\nabla u|^2}{2}$ which drives from a nonlinear optical model. The optical model described the propagation about self-trapped transverse magnetic field modes in a cylindrical optical fiber made from a self-focusing dielectric material. For more details, see \cite{Stuart12,Stuart13,Stuart C A,Zhou H S}.
\par
In recent years, there have been some studies on the particular quasilinear term $\phi$. For example, Jeanjean and R\u{a}dulescu \cite{Jeanjean L} extended the results in \cite{C.A. Stuart}. They studied the following second order quasilinear elliptic equation:
\begin{eqnarray}
\label{af3}
 \begin{cases}
 -\text{div}\left(\phi\left(\dfrac{u^2+|\nabla u|^2}{2}\right)\nabla u\right)+\phi\left(\dfrac{u^2+|\nabla u|^2}{2}\right)u
=f(u)+h,\;\;\;\;\;\; x\in \;\Omega ,\\
  u= 0,\;\;\;\;\;\; x\in \;\partial\Omega,
   \end{cases}
\end{eqnarray}
where $h\in L^2(\Omega)$ is non-negative, $\Omega$ is a bounded domain of $\mathbb{R}^N$ with $N\geq 1.$ They optimized the assumptions about the quasilinear part $\phi$. Specifically,  they assumed that $\phi$ satisfies the following conditions:\\
$(\phi_1)$ $\phi\in C([0,+\infty), \mathbb{R})$ and there exist some constants $0<\phi_{\min}\leq\phi_{\max}$ such that
$$\phi_{\min} \leq\phi(s)\leq\phi_{\max},\;\;\;\;\;s\in[0,+\infty);$$
$(\phi_2)$ the function $s\mapsto \Phi(s^2)$ is convex on $\mathbb{R}$, \\
where the condition $(\phi_2) $ ensures the ellipticity of (\ref{af3}) and is weaker than $(\Phi_2) $.
They obtained that  the problem $(\ref{af3})$ has a non-negative solution, when $h\gneqq 0$ and $f $ satisfies sublinear growth.  For $h\equiv0$, they got the result that the equation $(\ref{af3})$ has a  non-trivial solution if $f$  satisfies sublinear growth. Meanwhile, they  also claimed the nonexistence of solutions under some appropriate conditions. They deduced   that the problem $(\ref{af3})$  has at least one non-negative solution for $\|h\|_2$ sufficiently small and the solution is non-trivial if  $h \not\equiv0$, when $ f$ satisfies linear growth, and $\Phi$ satisfies $(\phi_2)$ and the following condition,\\
$(\phi_3)$ $\phi\in C([0,+\infty), \mathbb{R})$, there exist some constants $0<\phi_{\min}\leq\phi_{\max}$ such that
$$\phi_{\min} \leq\phi(s)\leq\phi_{\max},\;\;\;\;\;s\in[0,+\infty),$$
and there exists $\phi(\infty)>0$ such that $\phi(s)\rightarrow \phi(\infty)$ as $s\rightarrow \infty$.\\
They  also claimed that the problem $(\ref{af3})$ has  at least two non-trivial non-negative solutions for
 $h \not\equiv0$ and at least a non-trivial non-negative solution for $h \equiv0$ when $\phi$ satisfies $(\phi_3)$ and the following condition:\\
  $(\phi_4)$ the function $s\mapsto \Phi(s^2)$ is strictly convex on $\mathbb{R}$. \\
  What's more, combining the mountain pass theorem and some analytical skills, Pomponio and Watanabe \cite{Pomponio A 2021} got that  the equation $(\ref{af3})$ has a radial ground state solution $ u$ and $u\in C^{1,\sigma}$ for some $ \sigma\in (0,1)$, when  $h\equiv0$ and the nonlinearity satisfies a variant of Berestycki-Lions' conditions.
\par
In this paper, we first establish the existence of at least three non-trivial solutions for a class of quasilinear elliptic systems with Dirichlet boundary  value:
\begin{eqnarray}
\label{a}
 \begin{cases}
 -\text{div}\left(\phi_{1}\left(\dfrac{u^2+|\nabla u|^2}{2}\right)\nabla u\right)+\phi_{1}\left(\dfrac{u^2+|\nabla u|^2}{2}\right)u
 =\lambda G_{u}(x,u,v),\;\;\;\;\;\;x\in \;\Omega,\\
  -\text{div}\left(\phi_{2}\left(\dfrac{v^2+|\nabla v|^2}{2}\right)\nabla v\right)+\phi_{2}\left(\dfrac{v^2+|\nabla v|^2}{2}\right)v
 =\lambda G_{v}(x,u,v),\;\;\;\;\;\;x\in \;\Omega,\\
  u=v= 0,\;\;\;\;\;\;\;\;x\in \;\partial\Omega,\\
   \end{cases}
\end{eqnarray}
where $\Omega$ is a non-empty bounded open subset in $ \mathbb{R}^{N},$ $N\geq 1$, $|\Omega|:=\text{vol}(\Omega)=\int_\Omega 1dx<+\infty$,   $\lambda$ is a positive constant, $(u,v)\in W:=H_0^{1}(\Omega)\times H_0^{1}(\Omega)$ and $\Phi_i(s)=\int_0^s\phi_i(\tau)d\tau$, $i=1,2$. We assume that $\phi_i$  satisfies $(\phi_{4})$ and the following conditions:\\
$(\phi_{1})'$ $\phi_{i} \in C([0,+\infty),\mathbb{R})$  and there exist constants $0<\phi_{i}^0\leq\phi_{i}^m,$   such that
$$
\phi_i^{0}\leq\phi_{i}(s)\leq\phi_i^{m},\;\;\;\;\;\mbox{for all}\;s\in [0,+\infty).
$$
And the nonlinearity $G(x,\cdot,\cdot)$ is  $C^1$ in $\mathbb{R}\times \mathbb{R}$ for \text{a.e.} $x\in\Omega$ and  satisfies $G(\cdot,s,t)$ is measurable in $\Omega$ for any $(s,t)\in \mathbb{R}\times \mathbb{R}$ and the following  assumptions:\\
 $(G_{1})$ there exist non-negative constants $a_1,$ $a_2$ and $p,$ $q \in [1, 2^*)$ such that
 $$
 |G_s(x,s,t)|\leq a_1+a_2|s|^{p-1},\;\;\;\;\;\;\;\;\;\;|G_t(x,s,t)|\leq a_1+a_2|t|^{q-1},\;\;\;\;\;\;\;\forall (x, s, t)\in\Omega \times\mathbb{R}\times \mathbb{R};
 $$
 $(G_{2})$ $G(x,s,t)\geq 0 $ for every $ (x,s, t)\in\Omega \times\mathbb{R^+}\times \mathbb{R^+}$ and $G(x,0,0)=0$;\\
$(G_{3})$ there exist contants  $0<\alpha,$ $\beta<2$ and $b>0$ such that
 $$G(x,s,t)\leq b(1+|s|^\alpha+|t|^\beta),$$
 for almost every $x\in\Omega$ and  every $s,\;t\in\mathbb{R};$ \\
 $(G_{4})$ there exist  contants $\gamma_1>0$ and $\sigma_1>0$ with $\sigma_1>\gamma_1 \tau$, such that \\
 $$
\dfrac{ \inf_{x\in\Omega}G(x,\sigma_1,\sigma_1)}{\sigma_1^2}>2a_1\dfrac{A_1}{\gamma_1}+a_2A_2\gamma_1^{p-2}+a_2A_3\gamma_1^{q-2},
 $$
 where
  \begin{eqnarray}
\label{1a}
&&A_1=C_1\phi_m\left(\dfrac{1}{\phi_0}\right)^{1/2},\;\;\;
A_2=\dfrac{(C_p)^p\phi_m}{p}\left(\dfrac{2}{\phi_0}\right)^{p/2},\nonumber\\
&&
A_3=\dfrac{(C_q)^q\phi_m}{q}\left(\dfrac{2}{\phi_0}\right)^{q/2},\;\;\;
\tau=\dfrac{\sigma_1}{\left(\left[\Phi_1\left(\dfrac{\sigma_1^2}{2}\right)+\Phi_2\left(\dfrac{\sigma_1^2}{2}\right)\right]\dfrac{\pi^{N/2}}{\Gamma(1+N/2)}D^{N}\right)^{1/2}},
\end{eqnarray}
$a_1,$ $a_2$ are given in $(G_{1})$,  $\phi_0=\min\{\phi_1^0,\phi_2^0\}$, $\phi_m=\max\{\phi_1^m,\phi_2^m\}$, $D:=\sup_{x\in\Omega} \text{dist}(x,\partial \Omega)$,   $C_1$, $C_p$ and $C_q$ denote the embedding constants. Specifically, fixing $p\in [1,2^*]$, there exists the embedding $H_0^1(\Omega)\hookrightarrow L^p (\Omega)$ which is continuous and there exists a positive constant $C_p$ such that
 \begin{eqnarray}
\label{1}
\|u\|_p\leq C_p\|\nabla u\|_2.
\end{eqnarray}
In particular, the embedding $H_0^1(\Omega)\hookrightarrow L^p (\Omega)$ is compact when $p\in [1,2^*)$.  $H_0^{1}(\Omega)$ denotes the closure of $C_0^\infty$ in the Sobolev space $W^{1,2}(\Omega)$, equipped with the inner product and norm
$$
(u,v)=\int_{\Omega}\nabla u \cdot \nabla v dx,\;\;\;\;\;\;\| u\|_0=(u,u)^{1/2},
$$
where $W^{1,2}(\Omega)$ is defined as a space of all functions $u:\Omega\rightarrow\mathbb{R}$ satisfying $u\in L^2(\Omega)$ and the partial derivative $D^\alpha u\in L^2(\Omega)$ for $0\leq|\alpha|\leq1$.
$L^{p}(\Omega)$ $(1\leq p<\infty)$ denotes the Lebesgue space with the norm
$$\|u\|_{p}=\left(\int_{\Omega}|u|^{p}dx\right)^{1/p}.$$
 By simple caculation, there exists $x_0\in \Omega$ such that the ball $B(x_0,D)\subseteq \Omega$ and $\text{meas}(B(x_0,D))=\dfrac{\pi^{N/2}}{\Gamma(1+N/2)}D^{N}$.

\par
Presently, the existence of three  solutions  for  p-Laplacian equations and systems has been followed and many excellent results have been obtained. In \cite{Ricceri B2000}, Ricceri  studied a class of elliptic eigenvalue problems:
\begin{eqnarray}
\label{n2}
 \begin{cases}
 -\Delta u=\lambda [f(u)+\mu g(u)],\;\;\;\;\;\;x\in \;\Omega,\\
  u(x)= 0,\;\;\;\;\;\;\;\;x\in \;\partial\Omega,\\
   \end{cases}
\end{eqnarray}
 where $\Omega \subseteq \mathbb{R}^{N}$ is bounded, $N\geq 3$,  $\lambda$, $\mu \in \mathbb{R}$,  $f,\; g: \mathbb{R}\rightarrow \mathbb{R} $ are two continuous functions. Combining   \cite[Theorem 1]{Ricceri B1998} and  \cite[Corollary 1]{Pucci P1985}, he obtained   the existence of three non-trivial solutions for $(\ref{n2})$ when $\lambda>0$ and $|\mu|$ is small enough. \cite[Remarks 5.3 and 5.4]{Ricceri B2000} shows that  there is a trivial solution of  $(\ref{n2})$ for $g\equiv0$. That is,   there are only two non-trivial solutions for $g\equiv0$.  Then,  Ricceri improved the result of \cite[Theorem 3.1]{Ricceri B2000}  to  \cite [Theorem 1]{Ricceri B20002} which shows that the conclusion of  \cite[Theorem 3.1]{Ricceri B2000} holds  for each $\lambda$ in an open sub-interval  of an interval $I\subseteq \mathbb{R}$ by replacing  the convexity of the operator term with sequential weak lower semicontinuity. Under the additional assumption that the operator term is bounded on each bounded set of $\Omega$,  Ricceri \cite{Ricceri B2009} reached  the existence of three solutions for each $\lambda$ in a non-empty open interval $A$ and $A\subseteq I\subseteq \mathbb{R}$.  However, neither  \cite{Ricceri B20002} or   \cite{Ricceri B2009}  gave  the  further information on the size and location of the set  of $A$. Hereafter,  a new three critical points theorem was established in \cite{Ricceri B20092} ,  which ensures that  the conclusion of Theorem 1 in \cite{Ricceri B2009} actually holds for each compact interval A contained in $(a,+\infty)$, where $a$ is given, please refer to \cite{Ricceri B20092} for details.
Three critical points theorems for  non-differentiable functionals were established in \cite{Bonanno G 2008}, when the functions $f$, $g$ may be discontinuous on $u$, which extended the results previously in \cite{Ricceri B20002} from differentiable functionals to non-smooth functionals and determined a precise  interval of parameter $\lambda$.
In \cite[Theorem 3.5 and Theorem 3.6 ]{Bonanno G 2010}, Bonanno also got  a well-determined interval of parameter  $\lambda$  in the smooth and  the non-smooth framework for the functionals  which have at least three critical points under weaker regularity and compactness conditions.  In \cite{Bonanno Bisci 2011}, the existence of three non-trivial solutions for non-autonomous elliptic Dirichlet problems without any small perturbation of the nonlinearity is presented:
\begin{eqnarray}
\label{n3}
 \begin{cases}
 -\Delta u=\lambda f(x,u),\;\;\;\;\;\;x\in \;\Omega,\\
  u(x)= 0,\;\;\;\;\;\;\;\;x\in \;\partial\Omega,\\
   \end{cases}
\end{eqnarray}
where $f: \Omega \times\mathbb{R}\rightarrow \mathbb{R}$ is a function and satisfies the following assumptions:\\
$(F_{1})$ there exist two non-negative constants $a_1,$ $a_2$ and $q \in [1, 2^*)$ such that
 $$
 |f(x,t)|\leq a_1+a_2|t|^{q-1},\;\;\;\;\;\;\;\forall (x,  t)\in\Omega \times\mathbb{R};
 $$
 $(F_{2})$ $F(x,\xi)\geq 0 $ for every $ (x,\xi)\in\Omega \times\mathbb{R^+}$;\\
$(F_{3})$ there exist  contants  $0<s<2$ and $b>0$ such that
 $$F(x,\xi)\leq b(1+|\xi|^s),$$
 for almost every $x\in\Omega$ and every $\xi\in\mathbb{R};$ \\
 $(F_{4})$ there exist  contants $\gamma>0$ and $\delta>0$ with $\delta>\gamma \kappa$, such that \\
 $$
\dfrac{ \inf_{x\in\Omega}F(x,\delta)}{\delta^2}>a_1\dfrac{K_1}{\gamma}+a_2K_2\gamma^{q-2},
 $$
 where $a_1,$ $a_2$ are given in $(F_{1})$ and $\kappa$, $K_1$ and $K_2$ are given in  equation (10)  of \cite{Bonanno Bisci 2011}.
 \par
 They reached that the equation $(\ref{n3})$ has at least three non-trivial solutions in $H_0^1(\Omega)$ for each parameter $\lambda$ belonging to
 $$
  \Lambda_{\gamma, \delta }:=\left( \dfrac{2(2^N-1)}{D^2}\dfrac{\delta^2}{\inf_{x\in\Omega }F(x,\delta)}, \dfrac{2(2^N-1)}{D^2}\dfrac{1}{a_1K_1/\gamma+a_2K_2\gamma^{q-2}}\right).
  $$
The method they used is based on a recent three critical points theorem \cite[Theorem 3.6]{Bonanno G 2010} for differentiable functional, which can be rewrited as \cite[Theorem 2.1]{Bonanno Bisci 2011}.\\
 {\bf{Theorem I.}} (\cite [Theorem 2.1]{Bonanno Bisci 2011}) Let $X$ be a reflexive real Banach space, $\Phi: X\rightarrow \mathbb{R}$ be a coercive, continuously G\^{a}teaux differentiable and sequentially weakly lower semicontinuous functional whose G\^{a}teaux derivative admits a continuous inverse on $X^*$, $\Psi : X\rightarrow \mathbb{R}$ be a continuously G\^{a}teaux differentiable functional whose G\^{a}teaux derivative is compact such that
 $$
 \Phi(0)=\Psi(0)=0.
 $$
Assume that there exist $r>0$ and $\bar{x}\in X$, with $r<\Phi(\bar{x}),$ such that:\\
$(a_{1})$ $\dfrac{\sup_{\Phi(x)\leq r}\Psi(x)}{r}<\dfrac{\Psi(\bar{x})}{\Phi(\bar{x})}$; \\
$(a_{2})$ for each $\lambda \in \Lambda_\lambda :=\left( \dfrac{\Phi(\bar{x})}{\Psi(\bar{x})},\dfrac{r}{\sup_{\Phi(x)\leq r}\Psi(x)}\right)$, the functional $\Phi-\lambda \Psi$ is coercive.\\
 Then, for all $\lambda \in \Lambda_\lambda$ the functional $\Phi-\lambda \Psi$ has at least three distinct critical points in $X$.
\par
Moreover, using  variational methods and three non-trivial critical points theorem in \cite[Theorem 3.6]{Bonanno G 2010} or a more precise version of \cite[Theorem 2.1]{Bonanno Bisci 2011}, Bonanno-Heidarkhani-O'regan  proved  the existence of three non-zero  solutions for a gradient nonlinear Dirichlet elliptic system
 driven by a $(p,q)$-Laplacian operator in \cite{Bonanno G Heidarkhani 2011}, when $\Omega \subseteq \mathbb{R}^{N}$ is bounded and $p$, $q>N\geq1$.
 In \cite{Bonanno G 2017},  Averna and Bonanno  investigated  the existence of three non-trivial  weak solutions for the  Dirichlet elliptic system in \cite{Bonanno G Heidarkhani 2011},  when $N\geq 3$, $1<q\leq p <N$ and $f$ satisfies some appropriate assumptions.
  The results in \cite{Bonanno G 2017} determine an appropriate interval about the parameter $\lambda$ which improve the conclusions in \cite{Afrouzi G A}, \cite {Bonanno G Heidarkhani 2011} and \cite{Li C}.
 \par
 Meanwhile, the existence of three solutions for the Neumann problem is also  followed. In \cite{Bonanno G 2003 N}, applying Theorem 1 of  \cite{Ricceri B20002}, Bonanno and Candito  proved the existence of three non-trivial  weak solutions for the Neumann problem:
\begin{eqnarray}
\label{n5}
 \begin{cases}
 -\Delta_p u+a(x)|u|^{p-2}u
 =\lambda f(x,u),\;\;\;\;\;\;x\in \;\Omega,\\
  \dfrac{\partial u}{\partial\nu}= 0,\;\;\;\;\;\;\;\;x\in \;\partial\Omega,\\
   \end{cases}
\end{eqnarray}
where $\Omega \subseteq \mathbb{R}^{N}$ is bounded,   $f:\Omega\times \mathbb{R}\rightarrow\mathbb{R}$, $\lambda \in \mathbb{R_+}$, $p>N\geq1$, $a\in L^\infty(\Omega)$,  $ess\inf_{x\in\Omega}a(x)>0$ and $\nu$ is the outward unit normal to $\partial\Omega$. In \cite{D'Agu G},
applying the result which is obtained in \cite{Bonanno G 2010} and a more precise version of \cite[Theorem 3.2 ]{Bonanno G 2008},  D'Agui and Bisci   proved  the existence of at least three non-zero solutions for the problem of $(\ref{n5})$ and determined a precise interval of values of the parameter $\lambda$, when $\Omega \subseteq \mathbb{R}^{N}$ is bounded, $N\geq p>1$. Using \cite[Theorem 3.2 ]{Bonanno G 2008}, Bonanno- Bisci-R\u{a}dulescu  \cite{Bonanno G Bisci2011} studied the following   quasilinear elliptic Neumann problem:
\begin{eqnarray}
\label{n6}
 \begin{cases}
 -\text{div}(\alpha(|\nabla u|)\nabla u)+\alpha(|u|) u
 =\lambda f(x,u),\;\;\;\;\;\;x\in \;\Omega,\\
  \dfrac{\partial u}{\partial\nu}= 0,\;\;\;\;\;\;\;\;x\in \;\partial\Omega,\\
   \end{cases}
\end{eqnarray}
where $\Omega \subseteq \mathbb{R}^{N}$ is bounded, $N\geq 3$, and
$$
\varphi(t):=\left\{\begin{array}{l}
\alpha(|t|) t, \;\;\;\;\;\;\text {for } t\neq 0, \\
0,\;\;\;\;\;\;\text {for } t=0 \\
\end{array}\right.
$$
They obtained that the Neumann problem of $(\ref{n6})$ has at least three nontrivial solutions, when $ p_0>N\geq 3,$ where $$p_0:=\inf_{t>0}\dfrac{t\varphi(t)}{\Phi(t)},\;\;\;\;\Phi(t)=\int_0^t\varphi(s)ds.$$
  For a class of quasilinear elliptic systems involving the $p(x)$-Laplace operator with Neumann boundary condition, the existence of at least three  solutions is also established in \cite{Yin H Yang Z 2012} by using the three critical points theorem in \cite{Ricceri B2009}.
\par
Inspired by \cite{ Bonanno G Heidarkhani 2011, Bonanno Bisci 2011, Bonanno G Bisci2011,Yin H Yang Z 2012}, we investigate the following non-autonomous nonhomogeneous quasilinear elliptic Neumann problem:
\begin{eqnarray}
\label{aa}
 \begin{cases}
 -\text{div}\left(\phi_{1}\left(\dfrac{u^2+|\nabla u|^2}{2}\right)\nabla u\right)+\phi_{1}\left(\dfrac{u^2+|\nabla u|^2}{2}\right)u
 =\lambda G_{u}(x,u,v),\;\;\;\;\;\;x\in \;\Omega,\\
  -\text{div}\left(\phi_{2}\left(\dfrac{v^2+|\nabla v|^2}{2}\right)\nabla v\right)+\phi_{2}\left(\dfrac{v^2+|\nabla v|^2}{2}\right)v
 =\lambda G_{v}(x,u,v),\;\;\;\;\;\;x\in \;\Omega,\\
  \dfrac{\partial u}{\partial\nu}=\dfrac{\partial v}{\partial\nu}= 0,\;\;\;\;\;\;\;\;x\in \;\partial\Omega,\\
   \end{cases}
\end{eqnarray}
where $\Omega$ is a non-empty bounded open subset of the Euclidean space $ \mathbb{R}^{N},$ $N\geq 1$, $\lambda$ is a positive constant, $(u,v)\in X:=H^{1}(\Omega)\times H^{1}(\Omega)$.
We mainly study  the existence of at least three non-trivial solutions for  $(\ref{aa})$, under the assumptions $(\phi_1)'$, $(\phi_4)$, $(G_1)$-$(G_3)$ and the following condition:\\
$(G_{4})'$ there exist two positive contants $\gamma_2$ and $\sigma_2$, with $\sigma_2>\gamma_2 \bar{\tau}$ such that \\
 $$
\dfrac{ \inf_{x\in\Omega}G(x,\sigma_2,\sigma_2)}{\sigma_2^2}>2a_1\dfrac{\bar{A_1}}{\gamma_2}+a_2 \bar{A_2}\gamma_2^{p-2}+a_2\bar{A_3}\gamma_2^{q-2},
 $$
  \begin{eqnarray}
\label{2fa}
&&\bar{A_1}=c_1\phi_m\left(\dfrac{2}{\phi_0}\right)^{1/2},\;\;\;
\bar{A_2}=\dfrac{(c_p)^p\phi_m}{p}\left(\dfrac{2}{\phi_0}\right)^{p/2},\nonumber\\
&&
\bar{A_3}=\dfrac{(c_q)^q\phi_m}{q}\left(\dfrac{2}{\phi_0}\right)^{q/2},\;\;\;
\bar{\tau}=\tau=\dfrac{\sigma_2}{\left(\left[\Phi_1\left(\dfrac{\sigma_2^2}{2}\right)+\Phi_2\left(\dfrac{\sigma_2^2}{2}\right)\right]
\dfrac{\pi^{N/2}}{\Gamma(1+N/2)}D^{N}\right)^{1/2}},
\end{eqnarray}
 where $a_1,$ $a_2$ are given in $(G_{1})$, $\phi_0$, $\phi_m$ are given in $(\phi_1)'$, $c_1$, $c_p$ and $c_q$ are the embedding constants which are related with
 the embedding $H^1(\Omega)\hookrightarrow L^p(\Omega)$. That is, when $p\in [1,2^*)$, $H^1(\Omega)\hookrightarrow L^p(\Omega)$  and
the following  embedding inequality holds,
\begin{eqnarray}
\label{f7}
\|u\|_p\leq c_p\|u\|,
\end{eqnarray}
for every  $u\in H^{1}(\Omega)$. $H^{1}(\Omega)$ denotes the  space $W^{1,2}(\Omega)$ equipped with the norm
$$
\| u\|=\left[\int_\Omega(u^2+|\nabla u|^2)dx\right]^{1/2}.
$$
 \par
Next, we show our main results.
\par
\vskip2mm
 \noindent
 {\bf Theorem 1.1.} {\it Assume that $\phi_i$ $(i=1,2)$ and $G$ satisfy $(\phi_{1})'$, $(\phi_{4})$ and $(G_{1})$-$(G_{4})$ respectively. Then for $ N\geq 1$ and each parameter $\lambda$ belonging to
$$
\Lambda_{\gamma_1,\sigma_1}:=\left(\dfrac{\phi_m\sigma_1^2}{\inf_{x\in\Omega} G(x,\sigma_1,\sigma_1)},\dfrac{\phi_m}{\left( 2a_1\dfrac{A_1}{\gamma_1}+a_2 A_2{\gamma_1}^{p-2}+a_2 A_3{\gamma_1}^{q-2}\right)}\right),
$$
 the problem (\ref{a}) possesses at least three weak solutions in W, where $A_1$, $A_2$ and $A_3$ are given in $(\ref{1a})$.}
\par
Particularly, for $N=1$,  since the  embedding $H_0^1(\Omega)\hookrightarrow C^0(\overline{\Omega})$ is compact, it holds
\begin{eqnarray}
\label{aff7}
C:=\sup_{u\in H_0^1(\Omega)\backslash\{0\}}\dfrac{\|u\|_\infty}{\|u\|_0}<\infty.
\end{eqnarray}
Without the assumption $(G_1)$ and replacing $(G_4)$ with the following condition $(G_4)''$,\\
$(G_4)''$  there exist  contants $\gamma_1>0$ and $\sigma_1>0$ with $\sigma_1>\gamma_1 \tau$, such that \\
 $$
\dfrac{ \inf_{x\in\Omega}G(x,\sigma_1,\sigma_1)}{\sigma_1^2}>b(A_{11}\gamma_1^{-2}+A_{21}\gamma_1^{\alpha-2}+A_{31}\gamma_1^{\beta-2}),
 $$
 where
  \begin{eqnarray}
\label{1fa}
&&A_{11}=\phi_m|\Omega|,\;\;\;
A_{21}=\phi_m|\Omega|C^\alpha\left(\dfrac{2}{\phi_0}\right)^{\alpha/2},\nonumber\\
&&
A_{31}=b|\Omega|C^\beta \left(\dfrac{2}{\phi_0}\right)^{\beta/2},\;\;\;
\tau=\dfrac{\sigma_1}{\left(\left[\Phi_1\left(\dfrac{\sigma_1^2}{2}\right)+\Phi_2\left(\dfrac{\sigma_1^2}{2}\right)\right]|\Omega|\right)^{1/2}},
\end{eqnarray}
 we also get the existence of three non-zero solutions as the following.
\vskip2mm
 \noindent
 {\bf Theorem 1.2.} {\it Assume that $\phi_i$ $(i=1,2)$ and $G$ satisfy $(\phi_{1})'$, $(\phi_{4})$, $(G_{2})$-$(G_{3})$ and $(G_4)''$ respectively. Then for $N=1$ and each parameter $\lambda$ belonging to
$$
\bar{\Lambda}_{\gamma_1,\sigma_1}:=\left(\dfrac{\phi_m\sigma_1^2}{\inf_{x\in\Omega} G(x,\sigma_1,\sigma_1)},\dfrac{\phi_m}{b(A_{11}\gamma_1^{-2}+A_{21}\gamma_1^{\alpha-2}+A_{31}\gamma_1^{\beta-2})}\right),
$$
 the problem (\ref{a}) possesses at least three weak solutions in W, where $A_{11}$, $A_{21}$ and $A_{31}$ are given in $(\ref{1fa})$.}
\par
For  the problem (\ref{aa}), by replacing $(G_4)$ with $(G_4)'$, we can reach the following conclusion.
\vskip2mm
 \noindent
 {\bf Theorem 1.3.} {\it Assume that $\phi_i$ $(i=1,2)$ and $G$ satisfy $(\phi_{1})'$, $(\phi_{4})$ and $(G_{1})$-$(G_{3})$ and  $(G_{4})'$ respectively. Then for $N\geq 1$ and each parameter $\lambda$ belonging to
$$
\Lambda_{\gamma_2,\sigma_2}:=\left(\dfrac{\phi_m\sigma_2^2}{\inf_{x\in\Omega} G(x,\sigma_2,\sigma_2)},\dfrac{\phi_m}{\left( 2a_1\dfrac{\bar{A_1}}{\gamma_2}+a_2 \bar{A_2}{\gamma_2}^{p-2}+a_2 \bar{A_3}{\gamma_2}^{q-2}\right)}\right),
$$
 the problem (\ref{aa}) possesses at least three non-zero solutions in $X$, where $\bar{A_1}$, $\bar{A_2}$ and $\bar{A_3}$ are given in (\ref{2fa}).}
\par
Furthermore, for $N=1$,  the  embedding $H^1(\Omega)\hookrightarrow C^0(\overline{\Omega})$ is compact, and it holds
\begin{eqnarray}
\label{ff7}
c:=\sup_{u\in H^1(\Omega)\backslash\{0\}}\dfrac{\|u\|_{\infty}}{\|u\|}<\infty.
\end{eqnarray}
 We deduce the existence of three non-zero solutions for the problem $(\ref{aa})$ when $N=1$. The condition $(G_1)$ is not needed and we replace  $(G_4)'$ with the following condition $(G_4)'''$,\\
$(G_{4})'''$ there exist two positive contants $\gamma_3$ and $\sigma_3$, with $\sigma_3>\gamma_3 \hat{\tau}$ such that \\
 $$
\dfrac{ \inf_{x\in\Omega}G(x,\sigma_3,\sigma_3)}{\sigma_3^2}>b(\hat{A}_1\gamma_3^{-2}+\hat{A}_2\gamma_3^{\alpha-2}+\hat{A}_3\gamma_3^{\beta-2}),
 $$
 where
  \begin{eqnarray}
\label{2aa}
&&\hat{A}_1=\phi_m|\Omega|,\;\;\;
\hat{A}_2=\phi_m|\Omega|(c)^\alpha\left(\dfrac{2}{\phi_0}\right)^{\alpha/2},\nonumber\\
&&
\hat{A}_3=\phi_m|\Omega|(c)^\beta\left(\dfrac{2}{\phi_0}\right)^{\beta/2},\;\;\;
\hat{\tau}=\tau=\dfrac{\sigma_3}{\left(\left[\Phi_1\left(\dfrac{\sigma_3^2}{2}\right)+\Phi_2\left(\dfrac{\sigma_3^2}{2}\right)\right]|\Omega|\right)^{1/2}}.
\end{eqnarray}
and $a_1,$ $a_2$ are given in $(G_{1})$.
\vskip2mm
 \noindent
 {\bf Theorem 1.4.} {\it Assume that $\phi_i$ $(i=1,2)$ and $G$ satisfy $(\phi_{1})'$, $(\phi_{4})$ and $(G_{1})$-$(G_{3})$ and  $(G_{4})'''$ respectively. Then for $N= 1$ and each parameter $\lambda$ belonging to
$$
\Lambda_{\gamma_3,\sigma_3}:=\left(\dfrac{\phi_m\sigma_3^2}{\inf_{x\in\Omega} G(x,\sigma_3,\sigma_3)},\dfrac{\phi_m}{b( \hat{A}_1\gamma_3^{-2}+ \hat{A}_2{\gamma_3}^{\alpha-2}+ \hat{A}_3\gamma_3^{\beta-2})}\right),
$$
 the problem (\ref{aa}) possesses at least three weak solutions in $X$, where $\hat{A}_1$, $\hat{A}_2$ and $\hat{A}_3$ are given in (\ref{2aa}).}
\section{Preliminaries}
\label{}

Clearly, the corresponding functional of (\ref{a}) is defined on $W=H_0^{1}(\Omega)\times H_0^{1}(\Omega)$ by
\begin{eqnarray}
\label{a2}
J(u,v)=\int_{\Omega}\left[\Phi_{1}\left(\dfrac{u^2+|\nabla u|^2}{2}\right)+\Phi_{2}\left(\dfrac{v^2+|\nabla v|^2}{2}\right)\right]dx-\lambda\int_{\Omega}G(x,u,v)dx,
\end{eqnarray}
where the space  $W$ equiped with the norm
$$
\|(u,v)\|_{W}=\| u\|_0+\| v\|_0=\left(\int_{\Omega}|\nabla u|^2 dx\right)^{1/2}+\left(\int_{\Omega}|\nabla v|^2 dx\right)^{1/2}.
$$
With the similar arguements  in \cite[Appendix A.1]{Qi2023}, from $(\phi_1)'$ and $(G_1)$, one has $J\in C^1(W,\mathbb{R})$.    For any $(u,v)$ and $(u_{1},v_{1})\in W,$ we have
\begin{eqnarray}
\label{a4}
\langle J'(u,v),(u_{1},v_{1})\rangle&=&\int_{\Omega}\left[\phi_{1}\left(\dfrac{u^2+|\nabla u|^2}{2}\right)(uu_{1}+\nabla u\cdot \nabla u_{1})+\phi_{2}\left(\dfrac{v^2+|\nabla v|^2}{2}\right)(vv_{1}+\nabla v\cdot \nabla v_{1})\right]dx\nonumber\\
&&-\int_{\Omega}\left[G_{u}(x,u,v)u_{1}+G_{v}(x,u,v)v_{1}\right]dx\nonumber\\
&=&\langle J_{u}(u,v),u_{1}\rangle+\langle J_{v}(u,v),v_{1}\rangle.
\end{eqnarray}
Similarly, the corresponding functional of (\ref{aa}) is defined on $X=H^{1}(\Omega)\times H^{1}(\Omega)$ by
\begin{eqnarray*}
\label{aa2}
L(u,v)=\int_{\Omega}\left[\Phi_{1}\left(\dfrac{u^2+|\nabla u|^2}{2}\right)+\Phi_{2}\left(\dfrac{v^2+|\nabla v|^2}{2}\right)\right]dx-\lambda\int_{\Omega}G(x,u,v)dx.
\end{eqnarray*}
The space $X=H^{1}(\Omega)\times H^{1}(\Omega)$ equiped with the norm
$$
\|(u,v)\|_X=\|u\|+\|v\|=\left[\int_\Omega(u^2+|\nabla u|^2)dx\right]^{1/2}+\left[\int_\Omega (v^2+|\nabla v|^2)dx\right]^{1/2}.
$$
From $(\phi_1)'$ and $(G_1)$, it is easy to know  $L\in C^1(X,\mathbb{R})$. For any $(\bar{u},\bar{v})$
and $(u_{2},v_{2})\in X,$ we have
\begin{eqnarray}
\label{aa4}
\langle L'(\bar{u},\bar{v}),(u_{2},v_{2})\rangle&=&\int_{\Omega}\left[\phi_{1}\left(\dfrac{\bar{u}^2+|\nabla \bar{u}|^2}{2}\right)(\bar{u}u_{2}+\nabla \bar{u}\cdot \nabla u_{2})+\phi_{2}\left(\dfrac{\bar{v}^2+|\nabla \bar{v}|^2}{2}\right)(\bar{v}v_{2}+\nabla \bar{v}\cdot \nabla v_{2})\right]dx\nonumber\\
&&-\int_{\Omega}\left[G_{u}(x,\bar{u},\bar{v})u_{2}+G_{v}(x,\bar{u},\bar{v})v_{2}\right]dx\nonumber\\
&=&\langle L_{\bar{u}}(\bar{u},\bar{v}),u_{2}\rangle+\langle L_{\bar{v}}(\bar{u},\bar{v}),v_{2}\rangle
\end{eqnarray}
\par
First, we recall a convergence result \cite[Lemma 6]{Landes R 1980} about the monotone operators.
\vskip2mm
 \noindent
 {\bf Proposition 2.1.}
 {\it Let $X$ be a finite dimensional real Hilbert space with norm $| \cdot |$ and inner product $\langle \cdot, \cdot\rangle$. Assume that $\beta \in C(X, X)$ which is strictly monotone, that is
 \begin{eqnarray*}
\label{}
\langle \beta(\eta)- \beta(\eta_1),\eta-\eta_1\rangle>0,\;\;\;\;\mbox{for any}\;\eta, \eta_1\in X\;  \mbox{with } \eta\neq\eta_1.
\end{eqnarray*}
Setting $\{\eta_n\}\subset X$, $\eta\in X$ and
$$
\lim_{n\rightarrow +\infty}\langle \beta(\eta_n)- \beta(\eta),\eta_n-\eta\rangle=0.
$$
Then $\{\eta_n\}$ converges to $\eta$ in $X$.}
\par
 Similar to \cite[Lemma 3.4]{Jeanjean L}, we have the following conclusion,
 \vskip2mm
 \noindent
 {\bf Lemma 2.2.}
 {\it Assume that $(\phi_1)'$ and $(\phi_4)$ hold. Let $\beta : \mathbb{R}^{N+1}\rightarrow \mathbb{R}^{N+1}$ given by
 $$
\beta(\eta)=\phi\left(\dfrac{|\eta|^2}{2}\right)\eta.
 $$\\
 Then $ \beta$ is  strictly monotone.\\}
 {\bf Proof.} Let $ \eta=(u,\nabla u)$ and $ \eta_1=(\bar{u},\nabla \bar{u})$. Following from $(\phi_4)$, it is easy to see that
 $ \phi\left(\dfrac{t^2}{2}\right)t $ is strictly increasing on $\mathbb{R}$. In the view of \cite[Lemma 3.4]{Jeanjean L},  it holds that
  \begin{eqnarray}
\label{ab1}
&& \langle\phi\left(\dfrac{|\eta|^2}{2}\right)\eta-\phi\left(\dfrac{|\eta_1|^2}{2}\right)\eta_1,\eta-\eta_1\rangle>0, \;\;\; \text{for any}\;\eta,\eta_1\in \mathbb{R}^{N+1}.
\end{eqnarray}
  Then, we obtain that $\beta $  strictly monotone.
  \qed
   \vskip2mm
 \noindent
 {\bf Remark 2.3.} {\it In this paper, the quasilinear term depends on a function of $\dfrac{u^2+|\nabla u|^2}{2}$ but not $|\nabla u|$. This particular quasilinear term makes our problems are different from the problems about $p$-Laplacian\cite{Bonanno Bisci 2011,Li C} and $\Phi$-Laplacian\cite{Bonanno G Bisci2011,Shokooh  Graef  2020}. To obtain the  continuity of $(I_1')^{-1}$ and $(I_2')^{-1}$, it is necessary for our arguments  to use the theory of  monotone operators and the skill of adding one dimension to space. The definition of $(I_1')^{-1}$ and $(I_2')^{-1}$ will be given in  section 3 and section 4.}
\section{Dirichlet boundary Problem}
\label{}
In this section, we will give the proofs of Theorem 1.1 and Theorem 1.2. Firstly, we drive some lemmas those are useful for the proofs.
\par
In order to study problem (\ref{a}), we put the functional $I_1$, $\Upsilon_1 :H_0^{1}(\Omega)\times H_0^{1}(\Omega)\rightarrow \mathbb{R}$ and
\begin{eqnarray}
\label{a3}
I_1(u,v)=\int_{\Omega}\left[\Phi_{1}\left(\dfrac{u^2+|\nabla u|^2}{2}\right)+\Phi_{2}\left(\dfrac{v^2+|\nabla v|^2}{2}\right)\right]dx
\end{eqnarray}
and
\begin{eqnarray}
\label{a33}
\Upsilon_1(u,v)=\int_{\Omega}G(x,u,v)dx\;\;\;\;\;\forall \;(u,v)\in W.
\end{eqnarray}
\vskip2mm
 \noindent
 {\bf 3.1. Proofs of Theorem 1.1. }
 \vskip2mm
 \noindent
 {\bf Lemma 3.1.}
 {\it Assume that $(\phi_1)'$ and $(\phi_4)$ hold. Then functional $I_1\in C^1 (W, \mathbb{R}) $ is sequentially weakly lower semicontinuous, coercive and whose G\^{a}teaux derivative $I_1'$ admits a continuous inverse $(I_1')^{-1}$ on the dual space $W^*$ of $W.$}\\
{\bf Proof.} Similar to \cite[Lemma 3.2 (ii)]{Jeanjean L}, from $(\phi_4)$,  it is easy to know that
\begin{eqnarray}
\label{aaab1}
I_1(u,v)-I_1(\bar{u},\bar{v})-\langle I_1'(\bar{u},\bar{v}), (u,v)-(\bar{u},\bar{v})\rangle\geq0,\;\;\;\text{for any}\; (u,v),\;(\bar{u},\bar{v})\in W,
\end{eqnarray}
and
$I_1$ is sequentially weakly lower semicontinuous, that is

$$
I_1(u,v)\leq\liminf_{n\rightarrow\infty}I_1(u_n,v_n),\;\;\;\text{for}\;(u_n,v_n)\rightharpoonup (u,v)\; \text{in}\;W.
$$
 Following from $(\phi_1)'$, (\ref{a3}) and the Mean value theorems for definite integrals, there exist $\xi_1^x\in \left(0,\dfrac{u^2+|\nabla u|^2}{2}\right)$ and $\xi_2^x\in \left(0,\dfrac{v^2+|\nabla v|^2}{2}\right)$ for every $x\in\Omega$, such that
\begin{eqnarray*}
\label{b1}
I_1(u,v)&=&\int_{\Omega}\left(\int_0^{{(u^2+|\nabla u|^2)}/{2}}\phi_{1}(s)ds+\int_0^{{(v^2+|\nabla v|^2)}/{2}}\phi_{2}(t)dt\right)dx\nonumber\\
&=&\int_{\Omega}\left[\phi_{1}(\xi_1^x)\dfrac{u^2+|\nabla u|^2}{2}+\phi_{2}(\xi_2^x)\dfrac{v^2+|\nabla v|^2}{2}\right]dx\nonumber\\
&\geq&\dfrac{\phi_0}{2}(\|\nabla u\|_2^2+\|\nabla v\|_2^2)\nonumber\\
&\geq&\dfrac{\phi_0}{4}\|(u,v)\|_W^2.
\end{eqnarray*}
Obviously, $I_1(u,v)\rightarrow +\infty$ as $\|(u,v)\|_W\rightarrow +\infty$. Hence, $I_1$ is coercive.
\par
Next, we claim that  $(I_1')^{-1}\in C(W^*, W).$
From $(\phi_1)'$ and  (\ref{a4}), we have
\begin{eqnarray*}
\label{b3}
\|I_1'(u,v)\|_{W^*}&=&\sup_{(u_1,v_1)\in W/\{0\}}\dfrac{\langle I_1'(u,v),(u_1,v_1)\rangle}{\|(u_1,v_1)\|_W}\nonumber\\
&\geq&\dfrac{\langle I_1'(u,v),(u,v)\rangle}{\|(u,v)\|_W}\nonumber\\
&=&\dfrac{\displaystyle{\int_{\Omega}\left[\phi_{1}\left(\dfrac{u^2+|\nabla u|^2}{2}\right)(u^2+|\nabla u|^2)+\phi_{2}\left(\dfrac{v^2+|\nabla v|^2}{2}\right)(v^2+|\nabla v|^2)\right]dx}}{\|(u,v)\|_W}\nonumber\\
&\geq&\dfrac{\phi_0\| u\|_0^2+\phi_0\| v\|_0^2}{\|(u,v)\|_W}\nonumber\\
&=&\dfrac{\phi_0}{2}{\|(u,v)\|_W},
\end{eqnarray*}
for any $ (u,v)\in W$.
This shows that   $I_1'$ is coercive in $W$. Let $\eta=(u,\nabla u)$, $\eta_1=(\bar{u},\nabla\bar{u})$, $\xi=(v,\nabla v)$ and $\xi_1=(\bar{v},\nabla \bar{v})$.
Following from (\ref{a4}), (\ref{ab1}) and (\ref{a3}),
 it holds that
 \begin{eqnarray*}
\label{bbb3}
&&\langle I_1'(u,v)-I_1'(\bar{u},\bar{v}), (u,v)-(\bar{u},\bar{v})\rangle\nonumber\\
&=&\int_{\Omega}\phi_{1}\left(\dfrac{u^2+|\nabla u|^2}{2}\right)[u(u-\bar{u})+\nabla u\cdot (\nabla u-\nabla \bar{u})]dx\nonumber\\
&&+\int_{\Omega}\phi_{2}\left(\dfrac{v^2+|\nabla v|^2}{2}\right)[v(v-\bar{v})+\nabla v\cdot (\nabla v-\nabla \bar{v})]dx\nonumber\\
&&-\int_{\Omega}\phi_{1}\left(\dfrac{(\bar{u})^2+|\nabla \bar{u}|^2}{2}\right)[\bar{u}(u-\bar{u})+\nabla \bar{u}\cdot (\nabla u-\nabla \bar{u})]dx\nonumber\\
&&-\int_{\Omega}\phi_{2}\left(\dfrac{(\bar{v})^2+|\nabla \bar{v}|^2}{2}\right)[\bar{v}(v-\bar{v})+\nabla \bar{v}\cdot (\nabla v-\nabla \bar{v})]dx\nonumber\\
&=&\int_{\Omega}\left[\phi_{1}\left(\dfrac{|\eta|^2}{2}\right)\eta-\phi_{1}\left(\dfrac{|\eta_1|^2}{2}\right)\eta_1\right]\cdot
        (\eta-\eta_1)dx+\int_{\Omega}\left[\phi_{2}\left(\dfrac{|\xi|^2}{2}\right)\xi-\phi_{2}\left(\dfrac{|\xi_1|^2}{2}\right)\xi_1\right]\cdot (\xi-\xi_1)dx\nonumber\\
&>&0.
\end{eqnarray*}
 It holds that $I_1'$ is strictly monotone in $W$. From $(\phi_1)'$, it is easy to see that $I_1'$  is continuous. The continuity of $I_1'$ shows $I_1'$ is hemicontinuous. Together with \cite [Theorem 26]{Zeidler E 2013} and the fact that $I_1'$  is strictly monotone in $W$, it holds that there exists the inverse $(I_1')^{-1}$ of $I_1'$ and $(I_1')^{-1}$ is bounded in $W^*$. Next, we  claim that  $(I_1')^{-1}$  is continuous by proving  the inverse $(I_1')^{-1}$ of $I_1'$ is sequentially continuous. Put any sequence $\{w_n\}\subset W^*$  such that $w_n\rightarrow w \in W^*$. Let $(u_n,v_n)=(I_1')^{-1}(w_n)$ and $(u,v)=(I_1')^{-1}(w)$. We will prove that $(u_n,v_n)\rightarrow(u,v)$ in $W$. Following from $(I_1')^{-1}$ is bounded and $w_n\rightarrow w$ in $W^*$, it holds that $\{(u_n,v_n)\}$ is bounded in $W$. From the fact that $W$ is a reflexive Banach Space,  we have $(u_n,v_n)\rightharpoonup(u_0,v_0) $ in $W$, which shows that $u_n\rightharpoonup u_0$ in $H_0^1(\Omega)$ and $v_n\rightharpoonup v_0$ in $H_0^1(\Omega)$, respectively. From the boundness of $\{(u_n,v_n)\}$ and $w_n\rightarrow w$ in $W^*$, it holds that
$$
\langle w_n-w,(u_n,v_n)-(u_0,v_0)\rangle\rightarrow 0\;\;\;\mbox{as}\;n\rightarrow\infty,
$$
that is,
\begin{eqnarray}
\label{b4}
\langle I_1'(u_n,v_n)-I_1'(u,v),(u_n,v_n)-(u_0,v_0)\rangle\rightarrow 0\;\;\;\mbox{as}\;n\rightarrow\infty.
\end{eqnarray}
Since $(u_n,v_n)\rightharpoonup (u_0,v_0)$ in $W$, one has
\begin{eqnarray}
\label{b5}
\langle I_1'(u,v),(u_n,v_n)-(u_0,v_0)\rangle\rightarrow 0\;\;\;\mbox{as}\;n\rightarrow\infty
\end{eqnarray}
and
\begin{eqnarray}
\label{bb5}
\langle I_1'(u_0,v_0),(u_n,v_n)-(u_0,v_0)\rangle\rightarrow 0\;\;\;\mbox{as}\;n\rightarrow\infty.
\end{eqnarray}
Combining (\ref{b4}), (\ref{b5}) and (\ref{bb5}), we obtain that
\begin{eqnarray}
\label{b6}
0&=&\lim_{n\rightarrow\infty}[\langle I_1'(u_n,v_n),(u_n,v_n)-(u_0,v_0)\rangle-\langle I_1'(u_0,v_0),(u_n,v_n)-(u_0,v_0)\rangle]\nonumber\\
&=&\lim_{n\rightarrow\infty}\langle I_1'(u_n,v_n)-I_1'(u_0,v_0),(u_n,v_n)-(u_0,v_0)\rangle\nonumber\\
&=&\lim_{n\rightarrow\infty}\left\{\int_{\Omega}\left[\phi_{1}\left(\dfrac{u_n^2+|\nabla u_n|^2}{2}\right)u_n -\phi_{1}\left(\dfrac{u_0^2+|\nabla u_0|^2}{2}\right)u_0\right] (u_n-u_0)dx \right.\nonumber\\
&&+\int_{\Omega}\left[\phi_{1}\left(\dfrac{u_n^2+|\nabla u_n|^2}{2}\right)\nabla u_n -\phi_{1}\left(\dfrac{u_0^2+|\nabla u_0|^2}{2}\right)\nabla u_0\right]\cdot (\nabla u_n-\nabla u_0)dx \nonumber\\
&&+\int_{\Omega}\left[\phi_{2}\left(\dfrac{v_n^2+|\nabla v_n|^2}{2}\right)v_n -\phi_{2}\left(\dfrac{v_0^2+|\nabla v_0|^2}{2}\right)v_0\right] (v_n-v_0)dx \nonumber\\
&&\left.+\int_{\Omega}\left[\phi_{2}\left(\dfrac{v_n^2+|\nabla v_n|^2}{2}\right)\nabla v_n -\phi_{2}\left(\dfrac{v_0^2+|\nabla v_0|^2}{2}\right)\nabla v_0\right]\cdot (\nabla v_n-\nabla v_0)dx\right\}.
\end{eqnarray}
Define operators $\psi_i : H_0^1(\Omega)\rightarrow(H_0^1(\Omega))^*$, where $i=1,2$, and
\begin{eqnarray*}
\label{b7}
\langle \psi_1(u),u_1\rangle:=\int_{\Omega}\phi_{1}\left(\dfrac{u^2+|\nabla u|^2}{2}\right)(uu_1 +\nabla u \cdot \nabla u_1 )dx\;\;\;\forall u,\;u_1\in H_0^1(\Omega)
\end{eqnarray*}
and
\begin{eqnarray*}
\label{b8}
\langle \psi_2(v),v_1\rangle:=\int_{\Omega}\phi_{2}\left(\dfrac{v^2+|\nabla v|^2}{2}\right)(vv_1 +\nabla v \cdot \nabla v_1 )dx\;\;\;\forall v,\;v_1\in H_0^1(\Omega).
\end{eqnarray*}
From Lemma 2.2, it has that both  $\psi_1$ and $\psi_2$  are strictly monotone in $H_0^1(\Omega)$.
Together with  (\ref{b6}), we get
\begin{eqnarray}
\label{b9}
0&=&\lim_{n\rightarrow\infty}\left\{\int_{\Omega}\left[\phi_{1}\left(\dfrac{u_n^2+|\nabla u_n|^2}{2}\right)u_n -\phi_{1}\left(\dfrac{u_0^2+|\nabla u_0|^2}{2}\right)u_0\right] (u_n-u_0)dx \right.\nonumber\\
&&\left.+\int_{\Omega}\left[\phi_{1}\left(\dfrac{u_n^2+|\nabla u_n|^2}{2}\right)\nabla u_n -\phi_{1}\left(\dfrac{u_0^2+|\nabla u_0|^2}{2}\right)\nabla u_0\right]\cdot (\nabla u_n-\nabla u_0)dx\right\}
\end{eqnarray}
and
\begin{eqnarray}
\label{b10}
0&=&\lim_{n\rightarrow\infty}\left\{\int_{\Omega}\left[\phi_{2}\left(\dfrac{v_n^2+|\nabla v_n|^2}{2}\right)v_n -\phi_{2}\left(\dfrac{v_0^2+|\nabla v_0|^2}{2}\right)v_0\right] (v_n-v_0)dx \right.\nonumber\\
&&\left.+\int_{\Omega}\left[\phi_{2}\left(\dfrac{v_n^2+|\nabla v_n|^2}{2}\right)\nabla v_n -\phi_{2}\left(\dfrac{v_0^2+|\nabla v_0|^2}{2}\right)\nabla v_0\right]\cdot (\nabla v_n-\nabla v_0)dx\right\}.
\end{eqnarray}
From (\ref{bb5}), (\ref{b9}) and (\ref{b10}), it is easy to obtain that
\begin{eqnarray}
\label{b11}
0&=&\lim_{n\rightarrow\infty}\int_{\Omega}\left[\phi_{1}\left(\dfrac{u_n^2+|\nabla u_n|^2}{2}\right)u_n (u_n-u_0)
+\phi_{1}\left(\dfrac{u_n^2+|\nabla u_n|^2}{2}\right)\nabla u_n\cdot (\nabla u_n-\nabla u_0) \right]dx
\end{eqnarray}
and
\begin{eqnarray}
\label{b12}
0&=&\lim_{n\rightarrow\infty}\int_{\Omega}\left[\phi_{2}\left(\dfrac{v_n^2+|\nabla v_n|^2}{2}\right)v_n (v_n-v_0)
+\phi_{2}\left(\dfrac{v_n^2+|\nabla v_n|^2}{2}\right)\nabla v_n\cdot (\nabla v_n-\nabla v_0) \right]dx.
\end{eqnarray}
From \cite[Lemma 3.5]{Jeanjean L}, we get the conclusion that
\begin{eqnarray}
\label{bb12}
\nabla u_n(x)\rightarrow \nabla u(x),\;\;\;\text{a.e. in } \Omega.
\end{eqnarray}
In (\ref{aaab1}), let $u=u_0$, $v=v_0$, $\bar{u}=u_n $ and $\bar{v}=v_n$. Combining  (\ref{b11}) and (\ref{b12}),  it holds that
\begin{eqnarray}
\label{bbb44}
\limsup_{n\rightarrow \infty}I_1(u_n ,v_n  )\leq I_1(u_0,v_0).
\end{eqnarray}
Similar to \cite[Lemma 5.4]{Jeanjean L} or \cite[Lemma 2.3]{Pomponio A 2021}, let $j(s,t)=\Phi_1(s^2)+\Phi_2(t^2)$, $ p_n=h_n-h$  and $ q_n=g_n-g$, where
$$
 h_n=\left(\dfrac{u_n^2+|\nabla u_n|^2}{2}\right)^{1/2},\;\;\;\;\;h=\left(\dfrac{u_0^2+|\nabla u_0|^2}{2}\right)^{1/2}
 $$
 and
 $$
 g_n=\left(\dfrac{v_n^2+|\nabla v_n|^2}{2}\right)^{1/2},\;\;\;\;\;g=\left(\dfrac{v_0^2+|\nabla v_0|^2}{2}\right)^{1/2}.
 $$
 Following from $(\phi_1)'$ and $(\phi_4)$, it is easy to know that the function $j$ is continuous, strictly convex on $\mathbb{R}$ and $j(0,0)=0$. Using \cite[Lemma 5.4]{Jeanjean L} and (\ref{bb12}),  it holds that
 $$
 \int_\Omega[ j(h+p_n,g+q_n)-j(p_n,q_n)-j(h,g)]dx\rightarrow 0,\;\;\;\mbox{as}\; n\rightarrow\infty.
 $$
 That is,
 $$
 I_1(u_n,v_n)-I_1(u_n-u_0,v_n-v_0)-I_1(u_0,v_0)\rightarrow 0,\;\;\;\mbox{as}\; n\rightarrow\infty.
 $$
 Together with (\ref{bbb44}) and the fact that $I_1\geq 0$, we obtain
 \begin{eqnarray*}
\label{}
0 &\leq& \liminf_{n\rightarrow\infty}I_1(u_n-u_0,v_n-v_0)\nonumber\\
&\leq&\limsup_{n\rightarrow\infty}[I_1(u_n,v_n)-I_1(u_0,v_0)]\nonumber\\
&=&\limsup_{n\rightarrow\infty}I_1(u_n,v_n)-I_1(u_0,v_0)\nonumber\\
&\leq& 0.
 \end{eqnarray*}
This implies
\begin{eqnarray}
\label{q1}
\lim_{n\rightarrow\infty}I_1(u_n-u_0,v_n-v_0)=0.
 \end{eqnarray}
 From (\ref{a3}), it holds
\begin{eqnarray*}
\label{}
I_1(u,v)&\geq&\dfrac{\phi_0}{2} \int_\Omega(u^2+|\nabla u|^2+v^2+|\nabla v|^2)dx
\geq0,\;\;\;\;\mbox{for any}\; (u,v)\in W.
\end{eqnarray*}
Put $u=u_n-u_0$ and $v=v_n-v_0$. Following (\ref{q1}), it has
 \begin{eqnarray}
\label{q2}
 0&=&\lim_{n\rightarrow\infty}I_1(u_n-u_0,v_n-v_0)\nonumber\\
 &\geq&\dfrac{\phi_0}{2} \lim_{n\rightarrow\infty}\int_\Omega[(u_n-u_0)^2+|\nabla u_n-\nabla u_0|^2+(v_n-v_0)^2+|\nabla v_n-\nabla v_0|^2]dx\nonumber\\
 &\geq&0.
 \end{eqnarray}
 We get $(u_n,v_n)\rightarrow (u_0,v_0)$ in $W$. Together with  the continuity of $I_1'$, this shows  $I_1'(u_n,v_n)\rightarrow I_1'(u_0,v_0)=I_1'(u,v)$ in $W^*$,  when $n\rightarrow \infty$. Since $I_1'$ is  continuous and strictly monotone in $W^*$, it has $(u_0,v_0)=(u,v)$. Therefore, we conclude that  $(I_1')^{-1}$ is continuous.
\qed
\vskip2mm
 \noindent
 {\bf Lemma 3.2.}
 {\it Assume that $(G_1)$-$(G_4)$  hold. Then the functional $\Upsilon_1\in C^1 (W, \mathbb{R}) $  with compact derivative. Moreover,
 \begin{eqnarray}
\label{abb13}
\langle \Upsilon_1'(u,v),(u_1,v_1)\rangle=\int_\Omega \Upsilon_{1u}(u,v) u_1dx+\int_\Omega \Upsilon_{1v}(u,v) v_1dx
\end{eqnarray}
for all $(u_1,v_1)\in W.$}\\
{\bf Proof.}
Following from $(G_1)$, (\ref{1}) and (\ref{a33}),  one has
\begin{eqnarray*}
\label{bbb13}
\Upsilon_1(u,v)&\leq&\int_\Omega \left[ a_1(|u|+|v|) +a_2\dfrac{|u|^p}{p}+a_2\dfrac{|v|^q}{q}\right]dx\nonumber\\
&\leq& a_1C_1\|\nabla u\|_2+a_1C_1\|\nabla v\|_2+\dfrac{a_2(C_p)^p}{p}\|\nabla u\|_2^p+\dfrac{a_2(C_q)^q}{q}\|\nabla v\|_2^q.
\end{eqnarray*}
 Thus, $\Upsilon_1(u,v)$ is well defined in $W$. From $(G_1)$, it is easy to know that $\Upsilon_1\in C^1 (W, \mathbb{R}) $ and the equation $(\ref{abb13})$ holds.
\par
To reach the compactness of $\Upsilon_1'$, we put any sequence $\{(u_n,v_n)\}\subset W$ which is bounded. Following from the reflexivity of $W$ and the Sobolev embedding Theorem, there exists a subsequence still denoted by $\{(u_n,v_n)\}$, such that $(u_n,v_n)\rightharpoonup (u_0,v_0)\in W$, and $u_n\rightarrow u_0$ in $L^{p}(\Omega)$ and $v_n\rightarrow v_0$ in $L^{q}(\Omega)$, respectively. Similar to \cite{Wang L Zhang X Fang H 2017}, following from $(\ref{abb13})$, the H\"{o}lder's inequality and the Sobolev embedding inequality, we obtain
\begin{eqnarray}
\label{11b13}
&&|\Upsilon_1'(u_n,v_n)-\Upsilon_1'(u_0,v_0),(u_1,v_1)|\nonumber\\
&=&|\int_\Omega [G_u(x,u_n,v_n)-G_u(x,u_0,v_0)]u_1dx+\int_\Omega [G_v(x,u_n,v_n)-G_v(x,u_0,v_0)]v_1dx|\nonumber\\
&\leq&\bar{C}\left\{\left[\int_\Omega |G_u(x,u_n,v_n)-G_u(x,u_0,v_0)|^{p/(p-1)}dx\right]^{(p-1)/p}\right.\nonumber\\
&&\left.+\left[\int_\Omega| G_v(x,u_n,v_n)-G_v(x,u_0,v_0)|^{q/(q-1)}dx\right]^{(q-1)/q}\right\}\|(u_1,v_1) \|_W,
\end{eqnarray}
for any $(u_1,v_1)\in W,$ where $\bar{C}$ is a positive constant.  Meanwhile, by (\ref{11b13}) and the continuity of $G_u$ and $G_v$, one  has
\begin{eqnarray}
\label{12b13}
\int_\Omega |G_u(x,u_n,v_n)-G_u(x,u_0,v_0)|^{p/(p-1)}dx\rightarrow 0,\;\;\;\mbox{as}\;n\rightarrow\infty
\end{eqnarray}
and
\begin{eqnarray}
\label{199b13}
\int_\Omega|G_v(x,u_n,v_n)-G_v(x,u_0,v_0)|^{q/(q-1)}dx\rightarrow 0,\;\;\;\mbox{as}\;n\rightarrow\infty.
\end{eqnarray}
Combining $(\ref{11b13})$, $(\ref{12b13})$ and $(\ref{199b13})$, we obtain
\begin{eqnarray*}
\label{20b13}
&&\|\Upsilon_1'(u_n,v_n)-\Upsilon_1'(u_0,v_0)\|_*\nonumber\\
&=&\sup_{(\phi_1,\phi_2)\in W, \|(\phi_1,\phi_2)\|=1}|\langle\Upsilon_1'(u_n,v_n)-\Upsilon_1'(u_0,v_0),(\phi_1,\phi_2)\rangle|\nonumber\\
&=&\sup_{(\phi_1,\phi_2)\in W, \|(\phi_1,\phi_2)\|=1}|\int_\Omega[G_u(x,u_n,v_n)-G_u(x,u_0,v_0)]\phi_1dx
+\int_\Omega[G_v(x,u_n,v_n)-G_v(x,u_0,v_0)]\phi_2dx|\nonumber\\
&\leq&\sup_{(\phi_1,\phi_2)\in W, \|(\phi_1,\phi_2)\|=1}\left\{\int_\Omega|[G_u(x,u_n,v_n)-G_u(x,u_0,v_0)]\phi_1|dx
+\int_\Omega|[G_v(x,u_n,v_n)-G_v(x,u_0,v_0)]\phi_2|dx\right\}\nonumber\\
&\leq&\sup_{(\phi_1,\phi_2)\in W, \|(\phi_1,\phi_2)\|=1}\left\{\left[\int_\Omega |G_u(x,u_n,v_n)-G_u(x,u_0,v_0)|^{p/(p-1)}dx\right]^{(p-1)/p}\|\phi_1\|_{p}\right.\nonumber\\
&&\left.+\left[\int_\Omega| G_v(x,u_n,v_n)-G_v(x,u_0,v_0)^{q/(q-1)}dx\right]^{(q-1)/q}\|\phi_2\|_{q}\right\}\nonumber\\
&\leq&\sup_{(\phi_1,\phi_2)\in W, \|(\phi_1,\phi_2)\|=1}\left\{C_p\left[\int_\Omega |G_u(x,u_n,v_n)-G_u(x,u_0,v_0)|^{p/(p-1)}dx\right]^{(p-1)/p}\|\nabla \phi_1\|_{2}\right.\nonumber\\
&&\left.+C_q\left[\int_\Omega| G_v(x,u_n,v_n)-G_v(x,u_0,v_0)|^{q/(q-1)}dx\right]^{(q-1)/q}\|\nabla \phi_2\|_{2}\right\}\nonumber\\
&\leq&C\left\{\left[\int_\Omega |G_u(x,u_n,v_n)-G_u(x,u_0,v_0)|^{p/(p-1)}dx\right]^{(p-1)/p}\right.\nonumber\\
&&\left.+\left[\int_\Omega| G_v(x,u_n,v_n)-G_v(x,u_0,v_0)|^{q/(q-1)}dx\right]^{(q-1)/q}\right\}\|(\phi_1,\phi_2) \|_W\nonumber\\
&=&0.
\end{eqnarray*}
Then,  we reach the conclusion that $\Upsilon_1'$ is compact.
\qed
\par
Next, we will claim the assumptions $(a_1)$ and $(a_2)$ in Theorem I are satisfied. Following from $(G_1)$, one has
\begin{eqnarray}
\label{b13}
G(x,\xi,\zeta)\leq a_1(|\xi|+|\zeta|) +a_2\dfrac{|\xi|^p}{p}+a_2\dfrac{|\zeta|^q}{q},\;\;\;\mbox{for}\;\mbox{every}\;(x,\xi,\zeta)\in(\Omega\times \mathbb{R}\times \mathbb{R}).
\end{eqnarray}
Let $r_1\in(0,+\infty)$ and
\begin{eqnarray}
\label{b14}
\chi(r_1):=\dfrac{\sup_{(u,v)\in I_1^{-1}((-\infty,r_1])}\Upsilon_1(u,v)}{r_1}.
\end{eqnarray}
Following from (\ref{1}) and (\ref{b13}),one has
\begin{eqnarray}
\label{b15}
\Upsilon_1(u,v)&=&\int_\Omega G(x,u,v)dx\nonumber\\
&\leq& a_1\|u\|_1+a_1\|v\|_1+a_2\dfrac{\|u\|_p^p}{p}+a_2\dfrac{\|v\|_q^q}{q}\nonumber\\
&\leq& a_1C_1\|\nabla u\|_2+a_1C_1\|\nabla v\|_2+a_2(C_p)^p\dfrac{\|\nabla u\|_2^p}{p}+a_2(C_q)^q\dfrac{\|\nabla v\|_2^q}{q}.
\end{eqnarray}
  Following from  (\ref{b1}), for every $(u,v)\in W$ and $ I_1(u,v)\leq r_1$, we have
\begin{eqnarray}
\label{b16}
I_1(u,v)\geq \dfrac{\phi_0}{2}(\|\nabla u\|_2^2+\|\nabla v\|_2^2),
\end{eqnarray}
\begin{eqnarray}
\label{b17}
\dfrac{1}{2}(\|\nabla u\|_2+\|\nabla v\|_2)^2\leq\|\nabla u\|_2^2+\|\nabla v\|_2^2\leq \dfrac{2r_1}{\phi_0}
\end{eqnarray}
and
\begin{eqnarray*}
\label{}
\|\nabla u\|_2+\|\nabla v\|_2\leq 2\left(\dfrac{r_1}{\phi_0}\right)^{1/2}.
\end{eqnarray*}
Together with $(\ref{b15})$, it holds
\begin{eqnarray}
\label{bb17}
\sup_{(u,v)\in I_1^{-1}((-\infty,r_1])}\Upsilon_1(u,v)\leq 2a_1C_1\left(\dfrac{r_1}{\phi_0}\right)^{1/2}+\dfrac{a_2(C_p)^p}{p} \left(\dfrac{2r_1}{\phi_0}\right)^{p/2} +\dfrac{a_2(C_q)^q}{q} \left(\dfrac{2r_1}{\phi_0}\right)^{q/2} .
\end{eqnarray}
Following from (\ref{b14}) and (\ref{bb17}), it holds that
\begin{eqnarray}
\label{b18}
\chi(r_1) \leq 2a_1C_1 \left(\dfrac{1}{\phi_0}\right)^{1/2}r_1^{-1/2}+\dfrac{a_2(C_p)^p}{p} \left(\dfrac{2}{\phi_0}\right)^{p/2} r_1^{p/2-1}+\dfrac{a_2(C_q)^q}{q} \left(\dfrac{2}{\phi_0}\right)^{q/2}r_1^{q/2-1},\;\;\ r_1>0.
\end{eqnarray}
Put
$$
u_{\sigma_1}(x):=\left\{\begin{array}{l}
0, \;\;\;\;\;\;\text {if } x\in \Omega\backslash B(x_0,D), \\
\sigma_1,\;\;\;\;\;\;\text {if } x\in B(x_0,D) \\
\end{array}\right.
$$
and
$$
v_{\sigma_1}(x):=\left\{\begin{array}{l}
0, \;\;\;\;\;\;\text {if } x\in \Omega\backslash B(x_0,D), \\
\sigma_1,\;\;\;\;\;\;\text {if } x\in B(x_0,D). \\
\end{array}\right.
$$
Clearly, $(u_{\sigma_1},v_{\sigma_1})\in W$ and
\begin{eqnarray}
\label{b19}
I_1(u_{\sigma_1},v_{\sigma_1})&=&\int_\Omega\left[\Phi_{1}\left(\dfrac{u_{\sigma_1}^2+|\nabla u_{\sigma_1}^2|}{2}\right)+\Phi_{2}\left(\dfrac{v_{\sigma_1}^2+|\nabla v_{\sigma_1}^2|}{2}\right)\right]dx\nonumber\\
&=&\int\limits_{B(x_0,D)}\left[\Phi_{1}\left(\dfrac{\sigma_1^2}{2}\right)+\Phi_{2}\left(\dfrac{\sigma_1^2}{2}\right)\right]dx\nonumber\\
&=&\left[\Phi_{1}\left(\dfrac{\sigma_1^2}{2}\right)+\Phi_{2}\left(\dfrac{\sigma_1^2}{2}\right)\right] \text{meas}(B(x_0,D))\nonumber\\
&=&\left[\Phi_{1}\left(\dfrac{\sigma_1^2}{2}\right)
+\Phi_{2}\left(\dfrac{\sigma_1^2}{2}\right)\right]\dfrac{\pi^{N/2}}{\Gamma(1+N/2)}D^{N}.
\end{eqnarray}
From $(G_2)$, we have
\begin{eqnarray}
\label{b20}
\Upsilon_1(u_{\sigma_1},v_{\sigma_1})&=&\int\limits_\Omega G(x,u_{\sigma_1},v_{\sigma_1})dx\nonumber\\
&\geq&\int\limits_{B(x_0,D)}G(x,\sigma_1,\sigma_1)dx\nonumber\\
&\geq& \inf\limits_\Omega G(x,\sigma_1,\sigma_1)\dfrac{\pi^{N/2}}{\Gamma(1+N/2)}D^{N}.
\end{eqnarray}
Combining  (\ref{a33}), (\ref{b19}), (\ref{b20}) and $(\phi_1)'$, we get
\begin{eqnarray}
\label{b21}
\dfrac{\Upsilon_1(u_{\sigma_1},v_{\sigma_1})}{I_1(u_{\sigma_1},v_{\sigma_1})}&\geq&\dfrac{\inf\limits_\Omega G(x,\sigma_1,\sigma_1)}{\left[\Phi_{1}\left(\dfrac{\sigma_1^2}{2}\right)
+\Phi_{2}\left(\dfrac{\sigma_1^2}{2}\right)\right]}\nonumber\\
&\geq& \dfrac{1}{\phi_m}\dfrac{\inf\limits_\Omega G(x,\sigma_1,\sigma_1)}{\sigma_1^2}.
\end{eqnarray}
Note that $\tau=\dfrac{\sigma_1}{[I_1(u_{\sigma_1},v_{\sigma_1})]^{1/2}}$. Then  $ I_1(u_{\sigma_1},v_{\sigma_1})=\left(\dfrac{\sigma_1}{\tau}\right)^2$. Taking into account that $\sigma_1 >\gamma_1 \tau$,
 it holds $\gamma_1^2< I_1(u_{\sigma_1},v_{\sigma_1}).$
Together with (\ref{b14}), (\ref{b18}), (\ref{b21}) and $(G_4)$, we get
\begin{eqnarray*}
\label{b22}
\chi(\gamma_1^2):&=&\dfrac{\sup_{(u,v)\in I_1^{-1}((-\infty,\gamma_1^2])}\Upsilon_1(u,v)}{\gamma_1^2}\nonumber\\
&\leq&\dfrac{2a_1C_1}{\gamma_1} \left(\dfrac{1}{\phi_0}\right)^{1/2}+\dfrac{a_2(C_p)^p}{p} \left(\dfrac{2}{\phi_0}\right)^{p/2} \gamma_1^{p-2}+\dfrac{a_2(C_q)^q}{q} \left(\dfrac{2}{\phi_0}\right)^{q/2}\gamma_1^{q-2}\nonumber\\
&=&\dfrac{1}{\phi_m}\left( 2a_1\dfrac{A_1}{\gamma_1}+a_2 A_2{\gamma_1}^{p-2}+a_2 A_3{\gamma_1}^{q-2}\right)\nonumber\\
&<&\dfrac{1}{\phi_m}\dfrac{\inf\limits_\Omega G(x,\sigma_1,\sigma_1)}{\sigma_1^2}\nonumber\\
&\leq&\dfrac{\Upsilon_1(u_{\sigma_1},v_{\sigma_1})}{I_1(u_{\sigma_1},v_{\sigma_1})}.
\end{eqnarray*}
Therefore, the condition $(a_1)$ in Theorem I is satisfied.
\par
Moreover, for every $u$, $v\in H_0^1(\Omega)$, $0<\alpha<2$ and $0<\beta<2$,  it is easy to know that
$$|u|^\alpha\in L^{2/\alpha}(\Omega),\;\;\;\;\;|v|^\beta\in L^{2/\beta}(\Omega).$$
  Together with  the H\"{o}lder's inequality, it has
\begin{eqnarray}
\label{b23}
\int_\Omega|u|^\alpha dx\leq \|u\|^\alpha_2  |\Omega|^{(2-\alpha)/2},\;\;\;\forall u\in H_0^1(\Omega)
\end{eqnarray}
and
\begin{eqnarray}
\label{b24}
\int_\Omega|v|^\beta dx\leq \|v\|^\beta_2  |\Omega|^{(2-\beta)/2},\;\;\;\forall v\in H_0^1(\Omega).
\end{eqnarray}
Following from (\ref{1}), (\ref{a2}), (\ref{b23}), (\ref{b24}) and  $(G_3)$,  we obtain that
\begin{eqnarray*}
\label{b25}
J(u,v)&\geq& \dfrac{\phi_0}{2}(\|u\|^2_2+\|\nabla u\|^2_2+\|v\|^2_2+\|\nabla v\|^2_2)-\lambda b\int_\Omega(1+|u|^\alpha+|v|^\beta)dx\nonumber\\
&\geq&\dfrac{\phi_0}{2}(\|u\|^2_2+\|\nabla u\|^2_2+\|v\|^2_2+\|\nabla v\|^2_2)-\lambda b|\Omega|-\lambda b \|u\|^\alpha_2  |\Omega|^{(2-\alpha)/2}-\lambda b\|v\|^\beta_2  |\Omega|^{(2-\beta)/2}\nonumber\\
&\geq&\dfrac{\phi_0}{2}(\|\nabla u\|^2_2+\|\nabla v\|^2_2)-\lambda b |\Omega|-\lambda b (C_2)^\alpha|\Omega|^{(2-\alpha)/2}\|\nabla u\|^\alpha_2-\lambda b (C_2)^\beta|\Omega|^{(2-\beta)/2}\|\nabla v\|^\beta_2\nonumber\\
&\geq&\dfrac{\phi_0}{4}\|(u, v)\|^2_W-\lambda b |\Omega|-\lambda b (C_2)^\alpha|\Omega|^{(2-\alpha)/2}\|\nabla u\|^\alpha_2-\lambda b (C_2)^\beta|\Omega|^{(2-\beta)/2}\|\nabla v\|^\beta_2,
\end{eqnarray*}
for any $(u,v)\in W$.
Therefore, $J(u,v)$ is coercive  for every positive constant $\lambda$. Particularly, $J(u,v)$ is coercive
for every $\lambda\in\Lambda_{\gamma_1,\sigma_1}\subseteq\Bigg(\dfrac{I_1(u_{\sigma_1},v_{\sigma_1})}{\Upsilon_1(u_{\sigma_1},v_{\sigma_1})},
\dfrac{\gamma_1^2}{\sup_{(u,v)\in I_1^{-1}((-\infty,\gamma_1^2])}\Upsilon_1(u,v)}\Bigg).$ This shows that the condition $(a_2)$ in Theorem I is satisified. Therefore,  the conditions of Theorem I are satisfied. That is, we can reach the conclusion that the problem (\ref{a}) has at least three non-trivial weak solutions which are the critical point of the functional $J$ for any $\lambda\in\Lambda_{\gamma_1,\sigma_1}.$
\qed
\vskip2mm
 \noindent
 {\bf 3.2. Proofs of Theorem 1.2. }
 \vskip2mm
 \noindent
 \par
Since $H_0^1(\Omega)\hookrightarrow C^0(\overline{\Omega})$ is compact, when $\Omega\in \mathbb{R}^{N}$ and is a bounded domain with $N=1.$
It holds the Sobolev embedding inequality
\begin{eqnarray}
\label{az1}
\|u\|_\infty\leq C\|u\|_0<+\infty,
\end{eqnarray}
where $u\in H_0^{1}(\Omega)$ and $\|u\|_\infty:=\sup_{x\in \overline{\Omega}}|u(x)|.$
It is easy to prove that both Lemma 3.1 and Lemma 3.2 hold.
Following from $(G_3)$ and $(\ref{az1})$, one has
\begin{eqnarray}
\label{az2}
\Upsilon_1(u,v)&=&\int_\Omega G(x,u,v)dx\leq b\int_\Omega(1+|u|^\alpha+|v|^\beta)dx \nonumber\\
&\leq& b|\Omega|+b|\Omega|\|u\|^\alpha_\infty+b|\Omega|\|v\|^\beta_\infty\nonumber\\
&\leq& b|\Omega|+b(C)^\alpha|\Omega|\|u\|_0^\alpha+b(C)^\beta|\Omega|\|v\|_0^\beta.
\end{eqnarray}
Assume that $I_1(u,v)\leq r_1$. Following  from (\ref{b16}) and (\ref{b17}), it holds
\begin{eqnarray}
\label{azz4}
\|u\|_0\leq \left(\dfrac{2}{\phi_0}\right)^{1/2}r_1^{1/2},\;\;\;\;\|v\|_0\leq \left(\dfrac{2}{\phi_0}\right)^{1/2}r_1^{1/2},\;\;\;\;\|u\|_0+\|v\| _0 \leq2\left(\dfrac{r_1}{\phi_0}\right)^{1/2}.
\end{eqnarray}
Combining $(\ref{az2})$, (\ref{azz4}) and $(G_3)$, it follows that
\begin{eqnarray*}
\label{az5}
\Upsilon_1(u,v)&\leq&b|\Omega|+ b|\Omega|(C)^\alpha\left(\dfrac{2}{\phi_0}\right)^{\alpha/2}r_1^{\alpha/2}
+b(C)^\beta|\Omega|\left(\dfrac{2}{\phi_0}\right)^{\beta/2}r_1^{\beta/2}
\end{eqnarray*}
for $I_1(u,v)\leq r_1$. Thus, one has
\begin{eqnarray*}
\label{az6}
\chi_1(r_1)&=&\dfrac{\sup_{(u,v)\in I_1^{-1}((-\infty,r_1])}\Upsilon_1(u,v)}{r_1}\nonumber\\
&\leq&b|\Omega|r_1^{-1}+b|\Omega|(C)^\alpha\left(\dfrac{2}{\phi_0}\right)^{\alpha/2}r_1^{\alpha/2-1}
+b|\Omega|(C)^\beta\left(\dfrac{2}{\phi_0}\right)^{\beta/2}r_1^{\beta/2-1}.
\end{eqnarray*}
Put
$$
u_{\sigma_1}(x):=\sigma_1 \;\;\;\;\;\; v_{\sigma_1}(x):=\sigma_1,  \;\;\;\;\;\; \text{for}\; x\in \Omega.
$$
Clearly, $(u_{\sigma_1},v_{\sigma_1})\in X$.
This implies that
\begin{eqnarray*}
\label{az8}
I_1(u_{\sigma_1},v_{\sigma_1})&=&\int_\Omega\left[\Phi_{1}\left(\dfrac{u^2_{\sigma_1}+|\nabla u_{\sigma_1}|^2}{2}\right)+\Phi_{2}\left(\dfrac{v^2_{\sigma_1}+|\nabla v_{\sigma_1}|^2}{2}\right)\right]dx\nonumber\\
&=&\int_\Omega\left[\Phi_{1}\left(\dfrac{\sigma_1^2}{2}\right)+\Phi_{2}\left(\dfrac{\sigma_1^2}{2}\right)\right]dx\nonumber\\
&=&\left[\Phi_{1}\left(\dfrac{\sigma_1^2}{2}\right)+\Phi_{2}\left(\dfrac{\sigma_1^2}{2}\right)\right]|\Omega|
\end{eqnarray*}
and
\begin{eqnarray*}
\label{bb21}
\dfrac{\Upsilon_1(u_{\sigma_1},v_{\sigma_1})}{I_1(u_{\sigma_1},v_{\sigma_1})}&\geq&\dfrac{\inf_\Omega G(x,\sigma_1,\sigma_1)}{\left[\Phi_{1}\left(\dfrac{\sigma_1^2}{2}\right)
+\Phi_{2}\left(\dfrac{\sigma_1^2}{2}\right)\right]}\nonumber\\
&\geq& \dfrac{1}{\phi_m}\dfrac{\inf_\Omega G(x,\sigma_1,\sigma_1)}{\sigma_1^2}.
\end{eqnarray*}
 \par
Following from $(G_4)''$,
put $r_1=\gamma_1^2$,
\begin{eqnarray*}
\label{az11}
\chi(\gamma_1^2):&=&\dfrac{\sup_{(u,v)\in I_1^{-1}((-\infty,\gamma_1^2])}\Upsilon_1(u,v)}{\gamma_1^2}\nonumber\\
&\leq&b |\Omega|\gamma_1^{-2}+b|\Omega| (C)^\alpha\left(\dfrac{2}{\phi_0}\right)^{\alpha/2} \gamma_1^{\alpha-2}+b|\Omega|(C)^\beta \left(\dfrac{2}{\phi_0}\right)^{\beta/2}\gamma_1^{\beta-2}\nonumber\\
&=&\dfrac{1}{\phi_m}\left( bA_{11}\gamma_1^{-2}+b A_{21}{\gamma_1}^{\alpha-2}+b A_{31}{\gamma_1}^{\beta-2}\right)\nonumber\\
&<&\dfrac{1}{\phi_m}\dfrac{G(x,\sigma_1,\sigma_1)}{\sigma_1^2}\nonumber\\
&\leq&\dfrac{\Upsilon_1 (u_{\sigma_1},v_{\sigma_1})}{I_1 (u_{\sigma_1}v_{\sigma_1})}.
\end{eqnarray*}
This claims that  the assumption $(a_1)$ in Theorem I holds.
\par
Moreover,   one has
\begin{eqnarray*}
\label{az13}
\int_\Omega|u|^\alpha dx\leq\|u\|_\infty^\alpha|\Omega|, \;\;\;\forall u\in H_0^{1}(\Omega)
\end{eqnarray*}
and
\begin{eqnarray*}
\label{az14}
\int_\Omega|v|^\beta dx\leq\|v\|_\infty^\beta|\Omega|, \;\;\;\forall v\in H_0^{1}(\Omega).
\end{eqnarray*}
Combining with  $(\ref{az1})$, we have
\begin{eqnarray}
\label{az15}
\int_\Omega|u|^\alpha dx\leq\|u\|_\infty^\alpha|\Omega|\leq (C)^\alpha\|u\|_0^\alpha|\Omega|, \;\;\;\forall u\in H_0^{1}(\Omega),
\end{eqnarray}
and
\begin{eqnarray}
\label{azz15}
\int_\Omega|v|^\beta dx\leq\|v\|_\infty^\beta|\Omega|\leq (C)^\beta\|v\|_0^\beta|\Omega|, \;\;\;\forall v\in H_0^{1}(\Omega).
\end{eqnarray}
Following from $(\phi_1)'$, $(G_3)$,  (\ref{a2}) (\ref{az15}) and (\ref{azz15}), we obtain that
\begin{eqnarray*}
\label{f20}
J(u,v)&=&\int_{\Omega}\left[\Phi_{1}\left(\dfrac{u^2+|\nabla u|^2}{2}\right)+\Phi_{2}\left(\dfrac{v^2+|\nabla v|^2}{2}\right)\right]dx-\lambda\int_{\Omega}G(x,u,v)dx\nonumber\\
&\geq&\int_{\Omega}\left[\Phi_{1}\left(\dfrac{u^2+|\nabla u|^2}{2}\right)+\Phi_{2}\left(\dfrac{v^2+|\nabla v|^2}{2}\right)\right]dx\nonumber\\
&&-\int_{\Omega}\lambda bdx-\int_{\Omega}\lambda b |u|^\alpha dx-\int_{\Omega}\lambda b |v|^\beta dx\nonumber\\
&\geq&\dfrac{\phi_0}{2}\int_{\Omega}(u^2+|\nabla u|^2+v^2+|\nabla v|^2)dx-\lambda b|\Omega|\nonumber\\
&&-\lambda b\|u\|_\infty^\alpha|\Omega| -\lambda b\|v\|_\infty^\beta|\Omega|\nonumber\\
&\geq&\dfrac{\phi_0}{2}(\|u\|_0^2+\|v\|_0^2)-\lambda b|\Omega|-\lambda b(C)^\alpha\|u\|_0^\alpha|\Omega| -\lambda b(C)^\beta\|v\|_0^\beta|\Omega|,\;\;\;\forall (u,v)\in\;H_0^{1}(\Omega)\times H_0^{1}(\Omega).
\end{eqnarray*}
Since $\alpha$, $\beta<2,$ $J$ is a coercive functional for every positive parameter $\lambda$, and so for every
\begin{eqnarray*}
\label{f21}
\lambda\in\bar{\Lambda}_{\gamma_1,\sigma_1}=\left(\dfrac{\phi_m\sigma_1^2}{\inf_{x\in\Omega} G(x,\sigma_1,\sigma_1)},\dfrac{\phi_m}{b(A_{11}\gamma_1^{-2}+A_{21}\gamma_1^{\alpha-2}+A_{31}\gamma_1^{\beta-2})}\right).
\end{eqnarray*}
This shows that $(a_2)$ in Theorem I  holds. All these imply  $J$ has at least three non-trivial critical points for any $\lambda\in\bar{\Lambda}_{\gamma_1,\sigma_1}$. We complete the proofs of Theorem 1.2.\\
\qed
 \par
\vskip2mm
 \noindent
\section{Neuman boundary problem}
\vskip2mm
 \noindent
  \par
 In this section, we will give the proofs of Theorem 1.3 and Theorem 1.4  by replacing $(G_4)$  with $(G_4)'$ and $(G_4)'''$ respectively. This proofs are similar to section 3, here we mainly state the proofs related to $(X,\|\cdot\|)$ and the embedding inequality (\ref{f7}) or (\ref{ff7})
\par
In order to study problem (\ref{aa}), we put the functional $I_2$, $\Upsilon_2: X \rightarrow \mathbb{R}$ and
\begin{eqnarray}
\label{ff}
I_2(u,v)=\int_{\Omega}\left[\Phi_{1}\left(\dfrac{u^2+|\nabla u|^2}{2}\right)+\Phi_{2}\left(\dfrac{v^2+|\nabla v|^2}{2}\right)\right]dx
\end{eqnarray}
and
\begin{eqnarray}
\label{fff1}
\Upsilon_2(u,v)=\int_{\Omega}G(x,u,v)dx,\;\;\;\;\;\forall \;(u,v)\in X.
\end{eqnarray}
{\bf 4.1. Proofs of Theorem 1.3. }
\vskip2mm
 \noindent
 {\bf Lemma 4.1.}
 {\it Assume that $(\phi_1)'$ and $(\phi_4)$ hold. Then functional $I_2\in C^1 (X, \mathbb{R}) $ is sequentially weakly lower semicontinuous, coercive, bounded on each bounded subset of  $X$, and whose G\^{a}teaux derivative $I_2'$ admits a continuous inverse $(I_2')^{-1}$ on the dual space $X^*$ of $X.$}\\
{\bf Proof.} The arguements are similar to Lemma 3.1.
Following from $(\phi_4)$, it holds that
$$
I_2(u,v)\leq\liminf_{n\rightarrow\infty}I_2(u_n,v_n).
$$
From $(\ref{ff})$ and $(\phi_1)'$, one has
\begin{eqnarray}
\label{4f}
I_2(u,v)&\geq&\dfrac{\phi_0}{2}(\|u\|_2^2+\|\nabla u\|_2^2+\|v\|_2^2+\|\nabla v\|_2^2)\nonumber\\
&=&\dfrac{\phi_0}{2}(\|u\|^2+\|v\|^2)\nonumber\\
&\geq&\dfrac{\phi_0}{4}\|(u,v)\|_X^2.
\end{eqnarray}
This implies $I_2$ is coercive.
\par
Next, we deduce that $(I_2')^{-1}\in C(X^*,X)$. From $(\ref{aa4})$ and $(\phi_1)'$, we have
\begin{eqnarray*}
\label{6f}
\dfrac{\langle I_2'(u,v),(u,v)\rangle}{\|(u,v)\|_X}&=&\dfrac{\displaystyle{\int_{\Omega}\left[\phi_{1}\left(\dfrac{u^2+|\nabla u|^2}{2}\right)(u^2+|\nabla u|^2)+\phi_{2}\left(\dfrac{v^2+|\nabla v|^2}{2}\right)(v^2+|\nabla v|^2)\right]dx}}{\|(u,v)\|_X}\nonumber\\
&\geq&\dfrac{\phi_0\| u\|^2+\phi_0\|v\|^2}{\| (u,v)\|_X}\nonumber\\
&\geq&\dfrac{\phi_0}{2}\| (u,v)\|_X
\end{eqnarray*}
This shows that $I_2'$ is coercive in $X$. Similar to Lemma 3.1, we also can obtain $(I_2')^{-1}$ is continuous.
\qed
\par
By the similar arguements of Lemma 3.2, we can get the following result.
\vskip2mm
 \noindent
 {\bf Lemma 4.2.}
 {\it Assume that $(G_1)$-$(G_3)$ and $(G_4)'$ hold. Then the functional $\Upsilon_2\in C^1 (X, \mathbb{R}) $  with compact derivative. Moreover,
 \begin{eqnarray*}
\label{bb13}
\langle \Upsilon_2'(u,v),(u_1,v_1)\rangle=\int_\Omega \Upsilon'_{2u}(u,v) u_1dx+\int_\Omega \Upsilon'_{2v}(u,v) v_1dx
\end{eqnarray*}
for all $(u_1,v_1)\in X.$}
\par
Similar to (\ref{abb13}), we can obtain
 \begin{eqnarray*}
\label{}
\langle \Upsilon_2'(u,v),(u_1,v_1)\rangle=\int_\Omega \Upsilon'_{2u}(u,v) u_1dx+\int_\Omega \Upsilon'_{2v}(u,v) v_1dx
\end{eqnarray*}
for all $(u_1,v_1)\in X.$
Next, we will show the conditions $(a_1)$ and $(a_2)$ in Theorem I are satisfied. Let $r_2\in(0,+\infty)$ and
\begin{eqnarray}
\label{f1}
\chi(r_2):=\dfrac{\sup_{(u,v)\in I_2^{-1}((-\infty,r_2])}\Upsilon_2(u,v)}{r_2}
\end{eqnarray}
Assume that $I_2(u,v)\leq r_2$, together with (\ref{4f}), it has
\begin{eqnarray}
\label{f4}
\dfrac{\phi_0}{2}(\|u\|^2+\|v\|^2)\leq r_2.
\end{eqnarray}
Then, we have
\begin{eqnarray}
\label{f5}
\dfrac{1}{2}(\|u\|+\|v\|)^2\leq\|u\|^2+\|v\|^2\leq \dfrac{2}{\phi_0}r_2
\end{eqnarray}
and
\begin{eqnarray}
\label{f6}
\|u\|\leq \left(\dfrac{2}{\phi_0}\right)^{1/2}r_2^{1/2},\;\;\;\;\|v\|\leq \left(\dfrac{2}{\phi_0}\right)^{1/2}r_2^{1/2},\;\;\;\;\|u\|+\|v\| \leq2\left(\dfrac{r_2}{\phi_0}\right)^{1/2}.
\end{eqnarray}
Following from $(G_1)$, (\ref{f7}), (\ref{fff1}) and $(\ref{f4})$-$(\ref{f6})$, we have that
\begin{eqnarray}
\label{f8}
\Upsilon_2(u,v)&\leq&a_1\|u\|_1+a_1\|v\|_1+\dfrac{a_2}{p}\|u\|_p^p+\dfrac{a_2}{q}\|v\|_q^q\nonumber\\
&\leq&a_1c_1(\|u\|+\|v\|)+\dfrac{a_2}{p}(c_p)^p\|u\|^p+\dfrac{a_2}{q}(c_q)^q\|v\|^q\nonumber\\
&\leq&2a_1c_1\left(\dfrac{1}{\phi_0}\right)^{1/2}r_2^{1/2}+\dfrac{a_2}{p}(c_p)^p\left(\dfrac{2}{\phi_0}\right)^{p/2}r_2^{p/2}
+\dfrac{a_2}{q}(c_q)^q\left(\dfrac{2}{\phi_0}\right)^{q/2}r_2^{q/2},
\end{eqnarray}
for $I_2(u,v)\leq r_2$.  Combining $(\ref{f1})$ and $(\ref{f8})$, it holds that
\begin{eqnarray*}
\label{f9}
\chi_2(r_2)&=&\dfrac{\sup_{(u,v)\in I_2^{-1}((-\infty,r_2])}\Upsilon_2(u,v)}{r_2}\nonumber\\
&\leq&2a_1c_1\left(\dfrac{1}{\phi_0}\right)^{1/2}r_2^{-1/2}+\dfrac{a_2}{p}(c_p)^p\left(\dfrac{2}{\phi_0}\right)^{p/2}r_2^{p/2-1}
+\dfrac{a_2}{q}(c_q)^q\left(\dfrac{2}{\phi_0}\right)^{q/2}r_2^{q/2-1}.
\end{eqnarray*}
\par
Put
$$
u_{\sigma_2}(x):=\left\{\begin{array}{l}
0, \;\;\;\;\;\;\text {if } x\in \Omega\backslash B(x_0,D), \\
\sigma_2,\;\;\;\;\;\;\text {if } x\in B(x_0,D) \\
\end{array}\right.
$$
and
$$
v_{\sigma_2}(x):=\left\{\begin{array}{l}
0, \;\;\;\;\;\;\text {if } x\in \Omega\backslash B(x_0,D), \\
\sigma_2,\;\;\;\;\;\;\text {if } x\in B(x_0,D). \\
\end{array}\right.
$$
Note that $\bar{\tau}=\dfrac{\sigma_2}{[I_2(u_{\sigma_2},v_{\sigma_2})]^{1/2}}$. Then, $ I_2(u_{\sigma_2},v_{\sigma_2})=\left(\dfrac{\sigma_2}{\bar{\tau}}\right)^2$. Taking into account that $\sigma_2 >\gamma_2 \bar{\tau}$,
 it holds $\gamma_2^2< I_2(u_{\sigma_2},v_{\sigma_2}).$
 \par
Similar to the proofs of Theorem 1.1,  from $(G_4)$, we get
\begin{eqnarray*}
\label{b22}
\chi_2(\gamma_2^2):&=&\dfrac{\sup_{(u,v)\in I_1^{-1}((-\infty,\gamma_2^2])}\Upsilon_1(u,v)}{\gamma_2^2}
<\dfrac{\Upsilon_2(u_{\sigma_2},v_{\sigma_2})}{I_2(u_{\sigma_2},v_{\sigma_2})}.
\end{eqnarray*}
Therefore, the condition $(a_1)$ in Theorem I is satisified.
\par
Moreover,  from the H\"{o}lder inequality and $\alpha$, $\beta<2$, we obtain
\begin{eqnarray}
\label{f18}
\int_\Omega|u|^\alpha dx\leq\|u\|_2^\alpha |\Omega|^{(2-\alpha)/2}, \;\;\;\forall u\in H^{1}(\Omega).
\end{eqnarray}
Combining with  $(\ref{f7})$, we have
\begin{eqnarray}
\label{f19}
\int_\Omega|u|^\alpha dx\leq (c_2)^\alpha\|u\|^\alpha |\Omega|^{(2-\alpha)/2}, \;\;\;\forall u\in H^{1}(\Omega).
\end{eqnarray}
Similar to $(\ref{f18})$ and $(\ref{f19})$, it also has
\begin{eqnarray}
\label{ff19}
\int_\Omega|u|^\beta dx\leq (c_2)^\beta\|u\|^\beta |\Omega|^{(2-\beta)/2}, \;\;\;\forall u\in H^{1}(\Omega).
\end{eqnarray}
Following from $(\phi_1)'$, $(G_3)$,  (\ref{a2}) (\ref{f19}) and (\ref{ff19}),  we obtain that
\begin{eqnarray*}
\label{f20}
L(u,v)&=&\int_{\Omega}\left[\Phi_{1}\left(\dfrac{u^2+|\nabla u|^2}{2}\right)+\Phi_{2}\left(\dfrac{v^2+|\nabla v|^2}{2}\right)\right]dx-\lambda\int_{\Omega}G(x,u,v)dx\nonumber\\
&\geq&\int_{\Omega}\left[\Phi_{1}\left(\dfrac{u^2+|\nabla u|^2}{2}\right)+\Phi_{2}\left(\dfrac{v^2+|\nabla v|^2}{2}\right)\right]dx\nonumber\\
&&-\int_{\Omega}\lambda bdx-\int_{\Omega}\lambda b |u|^\alpha dx-\int_{\Omega}\lambda b |u|^\beta dx\nonumber\\
&\geq&\dfrac{\phi_0}{2}\int_{\Omega}[u^2+|\nabla u|^2+v^2+|\nabla v|^2]dx-\lambda b|\Omega|\nonumber\\
&&-\lambda b(c_2)^\alpha\|u\|^\alpha |\Omega|^{(2-\alpha)/2}-\lambda b(c_2)^\beta\|v\|^\beta |\Omega|^{(2-\beta)/2}\nonumber\\
&\geq&\dfrac{\phi_0}{4}\|(u,v)\|_X^2-\lambda b|\Omega|-\lambda b(c_2)^\alpha\|u\|^\alpha |\Omega|^{(2-\alpha)/2}-\lambda b(c_2)^\beta\|v\|^\beta |\Omega|^{(2-\beta)/2},
\end{eqnarray*}
for any $ (u,v)\in\;H^{1}(\Omega)\times H^{1}(\Omega).$
Hence, $L$ is a coercive functional for every positive parameter $\lambda$. So for every
\begin{eqnarray*}
\label{f21}
\lambda\in\Lambda_{\gamma_2,\sigma_2}=\left(\dfrac{\phi_m\sigma_2^2}{\inf_{x\in\Omega} G(x,\sigma_2,\sigma_2)},\dfrac{\phi_m}{\left( 2a_1\dfrac{\bar{A_1}}{\gamma_2}+a_2 \bar{A_2}{\gamma_2}^{p-2}+a_2 \bar{A_3}{\gamma_2}^{q-2}\right)}\right)
\end{eqnarray*}
the condition  $(a_2)$ in Theorem I holds. All these imply  the problem (\ref{aa}) has at least three critical points for every $\lambda\in\Lambda_{\gamma_2,\sigma_2}$. We reach the result of Theorem 1.3.
\qed
\par
\vskip2mm
 \noindent
{\bf 4.2. Proofs of Theorem 1.4.  }
\vskip2mm
 \noindent
 \par
The proofs  are similar to Theorem 1.2. We will omit them here.
\qed
\vskip2mm
 \noindent
\par
Similar to \cite[Example 3.1]{Bonanno Bisci 2011}, there are some examples about Theorem 1.1.
\par
\vskip2mm
 \noindent
{\bf Example. }
Let $\Omega$ be a non-empty bounded open subset of the Euclidean space $\mathbb{R}^N$. Put $p$, $q\in (1,2^*)$, $\alpha$, $\beta <2$ and
$$
w:=\max\{1, \tau, (A_1+A_2+A_3)^{1/(p-2)}p^{1/(p-2)}, (A_1+A_2+A_3)^{1/(q-2)}q^{1/(q-2)} \}.
$$
Let $r_1$ be a positive constant with $r_1>w$ and $G: \Omega\times\mathbb{R}\times\mathbb{R}\rightarrow\mathbb{R}$ and satisfies
$$
G_u(x,y,z):=\left\{\begin{array}{l}
1+|y|^{p-1},
\;\;\;\;\;\;\;\;\;\;\;\; y\leq r_1,\\
1+r_1^{p-\alpha}|y|^{\alpha-1},
\;\;\;\;\;\; y> r_1,
\end{array}\right.
$$
and
$$
G_v(x,y,z):=\left\{\begin{array}{l}
1+|z|^{q-1},
\;\;\;\;\;\;\;\;\;\;\;\; z\leq r_1,\\
1+r_1^{q-\beta}|z|^{\beta-1},
\;\;\;\;\;\; z> r_1,
\end{array}\right.
$$
Clearly, $(G_1)$ and $(G_2)$ hold. Moreover, for any $(s,t)\in \mathbb{R}\times\mathbb{R}$, it has
\begin{eqnarray*}
\label{g1}
G(x,s,t)&=&\int_0^{s}G_u(x,y,t)dy+\int_0^{t}G_v(x,0,z)dz\nonumber\\
&\leq&\left(r_1+\dfrac{r_1^\alpha}{r_1}\right)\left( 1+|\xi|^{\max\{1,\alpha\}}\right)+\left(r_1+\dfrac{r_1^\beta}{r_1}\right)\left( 1+|\eta|^{\max\{1,\beta\}}\right)\nonumber\\
&\leq&2\max\left\{\left(r_1+\dfrac{r_1^\alpha}{r_1}\right), \left(r_1+\dfrac{r_1^\beta}{r_1}\right)\right\}\left(1+|\xi|^{\max\{1,\alpha\}}+|\eta|^{\max\{1,\beta\}} \right).
\end{eqnarray*}
This shows that $(G_3)$ holds.
Since $r_1>w$, it holds that
\begin{eqnarray*}
\label{g2}
\dfrac{G(x,r_1,r_1)}{r_1^2}=\dfrac{r_1^{p-2}}{p}+\dfrac{1}{r_1}+\dfrac{r_1^{q-2}}{q}+\dfrac{1}{r_1}>A_1+A_2+A_3.
\end{eqnarray*}
\par
Particularly, let $\Omega$ be an open ball with radius one in $\mathbb{R}^N$, $p=q=3\in [1,2^*)$  and $\alpha=\beta=\dfrac{3}{2}<2$.  Consider the following systems,
\begin{eqnarray}
\label{ex}
 \begin{cases}
 -\text{div}\left(\left[1+\left(1+\dfrac{u^2+|\nabla u|^2}{2}\right)^{-1/2}\right]\nabla u\right)+\left[1+\left(1+\dfrac{u^2+|\nabla u|^2}{2}\right)^{-1/2}\right]u
 =\lambda G_{u}(x,u,v),\;\;\;\;\;\;x\in \;\Omega,\\
  -\text{div}\left(\left[2+\left(1+\dfrac{v^2+|\nabla v|^2}{2}\right)^{-1/3}\right]\nabla v\right)+\left[2+\left(1+\dfrac{v^2+|\nabla v|^2}{2}\right)^{-1/3}\right]v
 =\lambda G_{v}(x,u,v),\;\;\;\;\;\;x\in \;\Omega,\\
  u=v= 0,\;\;\;\;\;\;\;\;x\in \;\partial\Omega.\\
   \end{cases}
\end{eqnarray}
It is easy to see that $\phi_0=1$ and $\phi_m=3$. Following from the equations $(6)$ and $(7)$ in \cite{Bonanno Bisci 2011}, we obtain that
\begin{eqnarray*}
&&\label{}
(A_1+A_2+A_3)^{1/(p-2)}p^{1/(p-2)}\nonumber\\
&=& (A_1+A_2+A_3)^{1/(q-2)}q^{1/(q-2)}\nonumber\\
&=&\left[C_1\phi_m\left(\dfrac{2}{\phi_0}\right)^{1/2}+\dfrac{\phi_m(C_p)^p}{p}\left(\dfrac{2}{\phi_0}\right)^{3/2}
          +\dfrac{\phi_m(C_q)^q}{q}\left(\dfrac{2}{\phi_0}\right)^{3/2}\right]\nonumber\\
&=&3(A_1+A_2+A_3)\nonumber\\
&=&3\left[3\sqrt{2}C_1+2\sqrt[2]{2^3}(C_3)^3\right]\nonumber\\
&\leq&3\left\{ 3\sqrt{2}\times\dfrac{\pi}{4}\sqrt[4]{3}+4\sqrt{2}\times\left[ \left(\dfrac{\pi^2}{2}\right)^{1/4}\times (8\pi)^{-1/2}\times(6)^{1/4}\right]^3\right\}\nonumber\\
&\approx&14.9.
\end{eqnarray*}
Put $r_1=\sigma_1=20>w$, and  $G_u$, $G_v: \Omega\times\mathbb{R}\times\mathbb{R}\rightarrow\mathbb{R}$ denoted by
$$
G_u(x,y,z):=\left\{\begin{array}{l}
1+y^2,
\;\;\;\;\;\;\;\;\;\;\;\; y\leq 20,\\
1+20\sqrt{20y},
\;\;\;\;\;\; y> 20,
\end{array}\right.
$$
and
$$
G_v(x,y,z):=\left\{\begin{array}{l}
1+z^{2},
\;\;\;\;\;\;\;\;\;\;\;\; z\leq 20,\\
1+20\sqrt{20z},
\;\;\;\;\;\; z>20.
\end{array}\right.
$$
Thus, for any $\lambda\in\left(\dfrac{r_1^2}{G(x,r_1,r_1)},\dfrac{1}{A_1+A_2+A_3}\right)=\left(\dfrac{3r_1}{2r_1^2+6},\dfrac{1}{A_1+A_2+A_3}\right)\subseteq(0.07,4.9)$, the problem (\ref{ex}) has at least three non-trivial solution in $H_0^1(\Omega)\times H_0^1(\Omega).$
\par
The examples about Theorem 1.2, Theorem 1.3 and Theorem 1.4 can be obtained similarly.
\qed
\vskip4mm
 \noindent
{\bf Acknowledgement}\\
This work is supported by Yunnan Fundamental Research Projects of China (grant No: 202301AT070465) and  Xingdian Talent
Support Program for Young Talents of Yunnan Province in China.
\vskip2mm
 \noindent
 {\bf Conflict of interest}\\
On behalf of all authors, the  authors states that there is no conflict of interest.\\

\end{document}